\documentclass[10pt]{amsart}
\usepackage{amssymb,amsmath,amscd,amsfonts,amsthm,mathrsfs}
\usepackage{verbatim} 
\usepackage[dvips]{graphics}
\input xy
\xyoption{all}
\usepackage[all]{xy}  
\usepackage{color}

\usepackage[latin1]{inputenc}

\newcommand{\N}{\mathbb{N}}

\newcommand{\R}{\mathbb{R}}
\newcommand{\C}{\mathbb{C}}
\newcommand{\Ac}{\mathcal{A}}
\newcommand{\Bc}{\mathcal{B}}
\newcommand{\Cc}{\mathcal{C}}
\newcommand{\Dc}{\mathcal{D}}
\newcommand{\Ec}{\mathcal{E}}
\newcommand{\Hc}{\mathcal{H}}

\newcommand{\Id}{{\rm{Id}} }
\newcommand{\tr}{{\rm{tr}} }
\newcommand{\Sc}{\mathcal{S}}
\def\cchi{\raisebox{.45 ex}{$\chi$}}

\newcommand{\ran}{\mathrm{Ran}}  
\newcommand{\supp}{\mathrm{supp}}  
\providecommand{\father}[1]{\overleftarrow{#1}}

\newcommand{\Dr}{\mathscr{D}}
\newcommand{\Er}{\mathscr{E}}
\newcommand{\Fr}{\mathscr{F}}
\newcommand{\Gr}{\mathscr{G}}
\newcommand{\Hr}{\mathscr{H}}
\newcommand{\Kr}{\mathscr{K}}

\newcommand{\Nr}{\mathscr{N}}

\newcommand{\Qr}{\mathscr{Q}}
\newcommand{\Vr}{\mathscr{V}}

\def\build#1_#2^#3{\mathrel{\mathop{\kern 0pt#1}\limits_{#2}^{#3}}}

\newtheorem{theorem}{Theorem}[section]
\newtheorem{proposition}[theorem]{Proposition}
\newtheorem{lemma}[theorem]{Lemma}
\newtheorem{remark}[theorem]{Remark}
\newtheorem{corollary}[theorem]{Corollary}

\numberwithin{equation}{section}

\def \bone{\mathbf{1}}

\def \rmi{{\rm i}}

\title[Asymptotic eigenvalue distribution for
  discrete Laplacians]{Hardy inequality and asymptotic eigenvalue distribution for
  discrete Laplacians} 

\begin{document}
\author{Sylvain Gol\'enia}
\address{Institut de Math\'ematiques de Bordeaux, Universit\'e
Bordeaux $1$, $351$, cours de la Lib\'eration
\\$33405$ Talence cedex, France}
\email{sylvain.golenia@u-bordeaux1.fr}
\date{Version of \today}
\begin{abstract}
In this paper we study in detail some spectral properties of the
magnetic discrete Laplacian. We identify its form-domain,
characterize the absence of essential spectrum and provide the
asymptotic eigenvalue distribution.
\end{abstract}

\maketitle

\tableofcontents

\section{Introduction} 

The uncertainty principle is a central point in quantum physics. It can be
expressed by the following Hardy inequality:
\begin{eqnarray}\label{e:Hardy}
\left(\frac{n-2}{2}\right)^2\int_{\R^n} \left| \frac{1}{|x|}
f(x)\right|^2\, dx \leq \, \int_{\R^n}|\nabla f|^2\, dx = \langle
f, -\Delta_{\R^n} f\rangle,
\mbox{ where } n\geq 3,
\end{eqnarray}
and $f\in \Cc^\infty_c(\R^n)$. Roughly speaking, the
Laplacian controls some local singularities of a potential. In this paper, we
investigate which potentials a discrete Laplacian
is able to control. Obviously, since the value of a potential on a
vertex has to be finite, we will not focus on local
singularities. However, unlike in the continuous case, we will control
potentials that explode at infinity. 

We start with some definitions and fix our notation for graphs. We refer to
\cite{CdV, Chu, MW} for surveys on the matter.
Let $\Vr$ be a countable set. Let $\Er:=\Vr\times \Vr\rightarrow [0,\infty)$
and assume that 
\[\Er(x,y)=\Er(y,x), \quad \mbox{ for all } x,y\in \Vr.\]
We say that $G:=(\Er,\Vr)$ is an
unoriented weighted graph with \emph{vertices}~$\Vr$ and
\emph{weighted edges}~$\Er$. In the setting of electrical networks,
the weights correspond to the conductances.
We say that $x,y\in \Vr$ are \emph{neighbors} if $\Er(x,y)\neq 0$
and denote it by $x\sim y$. We say that there is a \emph{loop} in $x\in \Vr$
if $\Er(x,x)\neq 0$.  The set of \emph{neighbors} of $x\in \Er$
is denoted by 
\[\Nr_G(x):=\{y\in \Er,  x\sim y\}.\]
The \emph{degree} of $x\in V$ is by definition $|\Nr_G(x)|$,
the number of neighbors of~$x$.
A graph is \emph{locally finite} if $|\Nr_G(x)|$ is finite for all $x\in
V$. We also need a weight on the vertices 
\[m:\Vr\rightarrow
(0,\infty).\]
Finally, as we are dealing with magnetic fields, we fix a
phase 
\[\theta:\Vr\times \Vr\rightarrow [-\pi, \pi], \mbox{ such that }
\theta(x,y)= - \theta(y,x).\]
We set $\theta_{x,y}:=\theta(x,y)$.   A graph is \emph{connected}, 
if for all $x,y\in V$, there exists an $x$-$y$-\emph{path}, i.e.,
there is a finite sequence 
\[(x_1,\dotsc,x_{N+1})\in \Vr^{N+1} \mbox{ such that }
x_1=x, \, x_{N+1}=y \mbox{ and } x_n\sim x_{n+1}, 
\]
for all $n\in\{1,\dotsc,N\}$. The minimal possible $N$ is called the
(unweighted) \emph{distance} between  $x$ and $y$.

We recall that a graph~$G$ is \emph{simple} if
$\Er$ has values in $\{0,1\}$, $m=1$,  $\theta=0$, and has no loop. A
\emph{bi-partite} graph is a graph whose vertex set can be partitioned
into two subsets in such a way that no two points in the same subset
are neighbors. Trees are bi-partite graphs.

In the sequel, we shall always consider (magnetic) 
graphs $G=(\Vr, \Er, m, \theta)$, which are  locally
finite, connected and have no loop. We also fix $\omega\in \Vr$ and
denote by $|x|$ the distance between $x$ and $\omega$.

We now associate a certain Hilbert space and some operators on it to a given 
graph $G=(\Vr, \Er, m, \theta)$. Let $\ell^2(G, m^2):=\ell^2(\Vr, m^2;
\C)$ be the set of functions $f:\Vr\rightarrow \C$, such that
$\|f\|^2:=\sum_{x\in   \Vr} m(x)^2|f(x)|^2$ is finite. The associated
scalar product is given by $\langle f, g\rangle= \sum_{x\in \Vr}  
m^2(x) \overline{f(x)}g(x)$, for $f,g\in \ell^2(\Vr, m^2)$. 
We also denote by $\Cc_c(\Vr)$ the set of functions $f:\Vr\rightarrow
\C$, which have finite support. We define the quadratic form:
\begin{align}\label{e:quadra}
\Qr(f,f):= \Qr_{\Er, \theta}(f,f):=\frac{1}{2}\sum_{x,y\in \Vr}
\Er(x,y)|f(x)-e^{\rmi \theta_{x,y}}f(y)|^2\geq 0, \mbox{ for 
} f\in \Cc_c(V).
\end{align}
It is closable and there exists a unique self-adjoint operator
$\Delta_{\Er, \theta}$, such that  
\[\Qr_{\Er, \theta}(f,f)= \langle f,
\Delta_{\Er, \theta} f\rangle, \mbox{ for 
} f\in \Cc_c(V)
\]
and $\Dc(\Delta_{\Er, \theta}^{1/2})=
\Dc(\Qr_{\Er, \theta})$, where the latter is the completion of
$\Cc_c(\Vr)$ under $\|\cdot\|^2 + \Qr_{\Er,   \theta}(\cdot,
\cdot)$. This operator is the \emph{Friedrichs extension} associated
to the form $\Qr_{\Er, \theta}$ (see Section \ref{s:friedrichs} for
its construction). It acts as follows: 
\begin{align}\label{e:Deltam}
\Delta_{\Er, \theta} f(x) :=\displaystyle \frac{1}{m^2(x)}\sum_{y\in \Vr} \Er(x,y)
(f(x)- e^{\rmi \theta_{x,y}}f(y)), \mbox{ for
} f\in \Cc_c(\Vr).
\end{align}
When $m=1$, it is essentially self-adjoint
on $\Cc_c(V)$ (see Section \ref{s:esssa} for further discussion).   
If $G$ is simple, we shall simply write $\Delta_G$. There exist other
definitions for the discrete Laplacian, e.g., \cite{CdV, Chu, MW}, the one
we study here is sometimes called the ``physical Laplacian''. 

 In $\ell^2(G,
m^2)$, we define the \emph{weighted degree} by
\[d_G(x):=\frac{1}{m^2(x)} \sum_{y\in \Vr} 
\Er(x,y).\] 
Given a function $V:\Vr\to \C$, we
denote  by $V(Q)$ the operator of multiplication by $V$. It is
elementary that $\Dc(d_G^{1/2}(Q))\subset\Dc(\Delta^{1/2}_{\Er,
  \theta})$. Indeed, one has: 
\begin{align}\nonumber
\langle f, \Delta_{\Er, \theta} f\rangle&= \frac{1}{2} \sum_{x\in \Vr} \sum_{y\sim
  x} \Er(x,y)|f(x)-e^{\rmi \theta_{x,y}}f(y)|^2
\\
\label{e:majo}
&\leq \sum_{x\in \Vr} \sum_{y\sim 
  x} \Er(x,y)(|f(x)|^2+|f(y)|^2) = 2\langle f, d_G(Q) f\rangle,
\end{align}
for $f\in\Cc_c(\Vr)$. 
This inequality also gives a necessary condition for the absence of
essential spectrum for $\Delta_{\Er, \theta}$ (see Corollary
\ref{c:essDeltanec}). In Proposition \ref{p:bipart}, we also prove
that, in general, the constant $2$ cannot be improved. It
is also easy to see 
that $\Delta_{\Er, \theta}$ is bounded if and only if $d_G(Q)$ is
(see Proposition \ref{p:bounded}). 

In this paper we are interested in minorating the Laplacian with the
help of the weighted degree. In the case of non-magnetic Laplacians,
one standard approach is to use isoperimetric inequalities. The
classical version gives estimates on the bottom of the spectrum  . This is not
adapted to our situation. We rely on a modified version. 
We define the following
isoperimetric constant associated to (the weighted degree of) $G$ by 
\begin{align*}
\alpha(G):=\inf_{W\subset \Vr,\, \sharp W<\infty} \frac{\langle
  \bone_{W}, \Delta_{\varepsilon, 0} \bone_{W} \rangle }{\langle
  \bone_{W}, d_G(Q) \bone_{W} \rangle},
\end{align*}
where $\bone_X$ denotes the characteristic function of $X$. By
\cite[page 14, line -4]{KL2} (see also \cite{Do, DK, Kel} and
references therein), where $\alpha$ reads $\alpha_{b,c,n}$, we obtain: 
\begin{align}\label{e:iso}
\left(1- \sqrt{1- \alpha^2(G)}\right)\langle f, d_G(Q) f\rangle \leq \langle f,
\Delta_{\Er, 0}\, f\rangle 
\leq \left(1+ \sqrt{1- \alpha^2(G)}\right)\langle f, d_G(Q) f\rangle,
\end{align}
for all $f\in \Cc_c(\Vr)$. So, if $\alpha(G)>0$, \eqref{e:majo} is
improved and we have: $\Dc(\Delta_{\Er,
  0}^{1/2})=\Dc(d_G^{1/2}(Q))$. In particular, by  Proposition
\ref{p:compa}, one sees that the essential spectrum $\sigma_{\rm
  ess}(d_G(Q))$ is empty if and only if $\sigma_{\rm
  ess}(\Delta_{\Er, 0})$ is.  We refer to Theorem \ref{t:equadom} for
the equality of the domains and to Proposition \ref{p:pertutree2} for
the stability of the essential spectrum.     

We point out that a converse is also true. Namely, Proposition
\ref{p:negans} ensures that if there is $a>0$ so that $a\langle f,
\Delta_{\Er, 0}\, f\rangle \geq \langle f, d_G(Q)  f\rangle$, for all $f\in
\Cc_c(\Vr)$, then $\alpha(G)>0$. 

Assume that $\alpha(G)>0$. Supposing that $\sigma_{\rm
  ess}(\Delta_{\Er, 0})=\emptyset$ (or equivalently that
$\lim_{|x|\rightarrow \infty} d_G(x)=+\infty$), the inequality
\eqref{e:iso} and the min-max principle, see Proposition
\ref{p:compa}, provide the  bound 
\[  \left(1- \sqrt{1- \alpha^2(G)}\right)\leq  \liminf_{\lambda \to
  \infty}\frac{N_\lambda(\Delta_{\Er,  0})}{N_\lambda(d_G(Q))}
\leq \limsup_{\lambda \to  \infty}\frac{N_\lambda(\Delta_{\Er,
 0})}{N_\lambda(d_G(Q))}\leq \left(1+ \sqrt{1- \alpha^2(G)}\right),\] 
where 
\[N_\lambda(A):=\dim \ran\, \bone_{(-\infty, \lambda]}(A),
\]
 for a
self-adjoint operator $A$. This estimate has to be refined so as to
give the asymptotic of eigenvalues and to deal with magnetic fields.    
Moreover, \eqref{e:iso} is not stable by small perturbation for the
question of the equality of the form-domain. For instance, take a
simple graph $G_1$, such that $\alpha(G_1)>0$ and the (simple)
half-line graph $G_2$. Note that $\alpha(G_2)=0$. Now connect the
disjoint union of $G_1$ and $G_2$ by one edge to obtain a new graph
$G$. One sees easily  that $\alpha(G)=0$ and
$\Dc(\Delta_G^{1/2})=\Dc(d_G^{1/2}(Q))$. 
This is why we seek a minoration by $a d_G(Q)-b$ for some $a,b>0$.   

On the other hand, given $m_0:\Vr\to (0, \infty)$, we know that the
Laplacian acting in $\ell^2(\Vr, m^2)$ is unitarily equivalent to a
Schr\"odinger operator acting in  $\ell^2(\Vr, m^2_0)$ (see
Proposition \ref{p:uni}). This has already been noticed  before, e.g.,
\cite{CTT, HK}. By extracting some positivity, we obtain our analog of the
Hardy inequality:   

\begin{proposition}\label{p:hardy}
Let $G=(\Vr, \Er, m_0, \theta)$ be a locally finite graph. Given
$m:\Vr\rightarrow (0, \infty)$, one has  
\begin{align}\label{e:hardy}
\langle f, V_m(Q) f\rangle
\leq \langle f, \Delta_{\Er, \theta} f  \rangle, \mbox{ for } f\in
\Cc_c(\Vr),
\end{align}
where  
\begin{align}\label{e:Vm}
V_m(x):= d_G(x) - W_m(x), \mbox{ with }
W_m(x):=\frac{1}{m_0^2(x)}\sum_{y\in \Vr} \Er(x,y) 
  \frac{m(y)}{m(x)}\frac{m_0(x)}{m_0(y)}.
\end{align} 
Moreover, if $G$ is bi-partite, we  get:
\begin{align*}
\langle f, \Delta_{\Er, \theta} f  \rangle \leq \langle f,
(d_G(Q)+W_m(Q)) f\rangle, \mbox{ for } f\in
\Cc_c(\Vr).
\end{align*}
\end{proposition} 
Note that by choosing $m=m_0$, we recover that $\Delta_{\Er, \theta}\geq
0$. Moreover, $V_m$ is independent of the magnetic field. The
minoration is also different from the Kato's inequality of
\cite{DM}. We stress that the inequality \eqref{e:hardy} is in some
cases trivial, e.g., Proposition \ref{p:V=0}. One has to find a
favorable situation in order to exploit it. This is the case for some
perturbations of weighted trees. We present our main result:     

\begin{theorem}\label{t:tree}
Let $G_\circ=(\Vr, \Er_\circ, m, \theta)$ be a weighted
tree. Assume that there is $\varepsilon_0\in (0,1)$, so that  
\begin{align}\label{e:hypC0old}
C_0:=\sup_{x\in \Vr}\max_{y\in \Vr}\,
  d^{\varepsilon_0 -1}_{G_\circ}(x) \Er_\circ(x,y)m^{-2}(x)<\infty
  \mbox{ and } C_1:=\inf_{x\in \Vr} d_{G_\circ}(x)>0.
\end{align}
Let $G=(\Vr, \Er, m, \theta)$ be a perturbed graph and 
 $V:\Vr\rightarrow \R$ be a potential, satisfying:
\begin{align}\label{e:tree0}
|V(x)|+\Lambda(x) = o(1+ d_{G_\circ}(x)), \mbox{ as } |x|\to \infty, \mbox{
  where } \Lambda(x):=  \frac{1}{m^2(x)}\sum_{y\sim x} |\Er(x,y)-\Er_\circ(x,y)|. 
\end{align}
 Then, one has that:
\begin{enumerate}
\item The quadratic form associated to 
 $\Delta_{\Er,\theta}+V(Q)$ on $\Cc_c(\Vr)$ is bounded from below by some
constant $-C$. We denote by $\Hc_\Fr$ the associated Friedrichs extension. 
\item For all $\varepsilon>0 $, there is
  $c_\varepsilon\geq 0$, so that 
\begin{eqnarray}\label{e:tree1}
(1- \varepsilon)\langle f, d_{G_\circ}(Q) f\rangle - c_\varepsilon \|f\|^2\leq
\langle f, \Hc_\Fr  f\rangle \leq
(1+ \varepsilon) \langle f, d_{G_\circ}(Q) f\rangle+ c_\varepsilon
\|f\|^2,
\end{eqnarray}
for $f\in  \Cc_c(\Vr)$. We have $\Dc(|\Hc_\Fr|^{1/2})=
\Dc((d_{G_\circ}(Q)^{1/2})$.  
\item The essential spectrum of $\Hc_\Fr$ is equal to 
  that of $\Delta_{\Er_\circ, \theta_\circ}$. 
\item The essential spectrum of
$\Hc_\Fr$ is empty if and only if $\lim_{|x|\rightarrow
  \infty} d_{G_\circ}(x)= +\infty$. In this case we 
obtain:
\begin{align}\label{e:asymp}
\lim_{N \to
  \infty}\frac{\lambda_N(\Hc_\Fr)}{\lambda_N(d_{G_\circ}(Q))}=1,
\end{align}
where $\lambda_N$ denotes the $N$-th eigenvalue counted with multiplicity.
\end{enumerate}
\end{theorem} 

The theorem will be proved in Section \ref{s:hardytree}. In \eqref{e:tree0}, we  have use the Landau's notation for the small $o$. Namely, $f(x)=o(g(x))$, as $|x|\to \infty$ if $f/g(x)$ tends to
$0$ as $|x|\to \infty$. Notice that the convergence given by $|x|\to
\infty$ corresponds that given by the filter generated by the complements of
finite sets and is independent of the choice of $\omega$. 

Note that we improve on the bound \eqref{e:majo}. We point out that Hypothesis
\eqref{e:hypC0old} is fulfilled by simple trees. We improve this
condition in \eqref{e:hypC0} and discuss it Remark \ref{r:hypC0}. 
Since $G_0$ is a
tree, we recall that 
$\Delta_{\Er_\circ, \theta}$ is unitarily equivalent to
$\Delta_{\Er_\circ, 0}$. However, $G$ is a priori not a tree (recall
that Zorn's Lemma ensures that every simple graph has a maximal 
subtree). Therefore it is interesting to observe that there is no
hypothesis on $\theta$.  We indicate that the inequality
\eqref{e:tree1} is valid for a larger class of perturbations (see
Proposition \ref{p:pertutree}).  

We point out that the first part of d), namely the absence of the
essential spectrum, has been studied in many works, e.g., \cite{Kel,
  KL2, KL, KLW}. They generalize some ideas of \cite{DK,
  Fuj}. Their approach is based on some isoperimetric estimates and
on the Persson's Lemma. The latter characterizes the infimum of the
essential spectrum.     

The asymptotic of eigenvalues is a novelty and was not considered in the
literature before. Here one should keep in mind that our approach is
different  from the one used in the continuous setting. Whereas one
usually relies 
on the Dirichlet-Neumann bracketing technique, by cutting the space
into boxes, it is hard to believe that such an approach would be
efficient here. Indeed, cutting the graph gives a perturbation which
is of the same size as the operator.  

We stress that one can prescribe any asymptotic of
eigenvalues by choosing a proper tree $G$ (and in fact $d_G$). We
mention that the spectral asymptotic estimates obtained in \cite{DM}
are for some operators with non-empty essential spectrum. They study
graphs which are equipped with a free action of a discrete group and
establish a bound on the  $\tr\, e^{-t\Delta_{\Er,  \theta}}$, where
the trace is adapted to a  fundamental domain.   

We turn to the question of the form-domain. We stress that we do not
suppose that the isoperimetric constant is non-zero.
To our knowledge, this is the first time that the
form-domain of the unbounded discrete Laplacian on a simple tree is
identified.  It is remarkable that the form-domain coincides with
that of $d_G(Q)$, a multiplication  operator. A useful consequence is
the stability of the essential spectrum, obtained in c). This is also
new. On the other hand, we stress that there are simple
bi-partite graphs, such that the  form-domain of the Laplacian is
different from that of $d_G(Q)$ (see Proposition   \ref{p:bipart}). In
this case, \eqref{e:tree1} is not fulfilled.    

Having the same form-domain does not necessarily ensure that the
domains are also equal.  In Proposition \ref{p:treedom}, we construct a
simple tree 
which is such an example. However, under some further hypotheses on
the graph, Theorem \ref{t:equadom} ensures that the domain of the
magnetic Laplacian is equal to that of $d_G(Q)$. In Proposition
\ref{p:treenotdom}, we give an example of a simple tree $T$, which has $0$
as associated isoperimetric constant and such that the domain of the
Laplacian is the same as that of the weighted degree. Moreover,
one obtains that $\sigma(\Delta_T)=\sigma_{\rm ac}(\Delta_T)=[0,\infty)$.

Finally we present the organization of this paper. In section
\ref{s:esssa} we provide a new criterion of essential
self-adjointness. Next, in section \ref{s:min-max}, we recall
some well-known facts about the min-max principle, its relation to
the bottom of the essential spectrum and compactness. Then, in Section
\ref{s:hardy}  we prove the Hardy inequality and discuss its
triviality. In Section \ref{s:hardytree}, we prove Theorem
\ref{t:tree} in the context of trees.  Next, in Section \ref{s:dom}
we discuss the question of the domain of the Laplacian on a general
graph and that of form-domain on bi-partite
graphs in Section \ref{s:dombip}. Perturbation theory is developed in Section
\ref{s:pertu}. Finally we provide two appendices, one
concerning the $\Cc^1$ regularity and another one concerning the 
Helffer-Sj\"{o}strand's formula.

\vskip1mm
\noindent\textbf{Notation:} We denote by $\N$ the non-negative
integers. In particular, $0\in \N$.  We set $\langle
x\rangle:=(1+x^2)^{1/2}$. Given a set~$X$ and $Y\subseteq X$ let
$\bone_Y\colon X\to\{0,1\}$ be the characteristic function of~$Y$. We
denote also by $Y^c$ the complement set of~$Y$ in~$X$. We consider
only separable complex Hilbert space. We denote by $\Bc(\Hr, \Kr)$,
the space of bounded operators between the Hilbert spaces $\Hr$ and
$\Kr$.

\vskip1mm
\noindent\textbf{Acknowledgments:} I thank heartily
Michel Bonnefont, Thierry Jecko, Matthias Keller, Daniel Lenz, Ognjen
Milatovic, Sergiu Moroianu, Elizabeth Strouse,  Fran\c coise Truc, and the
anonymous referee for fruitful
discussions and comments on the script.  

\section{General properties}
\subsection{A few words about the Friedrichs
  extension}\label{s:friedrichs} 
Given a dense subspace $\Dr$ of a Hilbert space $\Hr$ and a
non-negative symmetric operator $H$ on $\Dr$. We define  
the quadratic
from $\Qr(f,g):= \langle  f, H g \rangle + \langle f,g\rangle $ on $\Dr\times
\Dr$. Let $\Hr_1$ be the
completion of $\Dr$ under the norm associated to $\Qr$, namely by the
norm $\|\cdot\|_{\Qr}$ given by $\|\varphi\|^2_\Qr:=\Qr(\varphi)^2=\langle 
H\varphi, \varphi\rangle+\|\varphi\|^2$.  The domain of the
Friedrichs extension of $H$ is given by 
\begin{align*}
\Dc(H_\Fr)&=\{f\in \Hr_1, 
\Dr\ni g\mapsto \langle Hg,f\rangle+\langle g,f\rangle \mbox{ extends
  to a norm continuous function on } \Hr\}
\\
&= \Hr_1\cap \Dc(H^*).
\end{align*}
For each
$f\in\Dc(H_\Fr)$, there is a unique $u_f$ such that $\langle
Hg,f\rangle + \langle g,f\rangle  = \langle g,u_f\rangle$, by Riesz'
Theorem. The Friedrichs extension of $H$, is given by $H_\Fr
f=u_f-f$. It is a self-adjoint extension of $H$, e.g., \cite[Theorem
X.23]{RS}. Moreover $\Dc((H_\Fr)^{1/2})= \Hr_1$. In the sequel we drop
the notation with $\Fr$ when we refer to the Friedrichs extension of
the Laplacian, i.e., $\Delta_{\Er, \theta}=(\Delta_{\Er, \theta})_\Fr$.

It remains to describe the domain of the adjoint of a discrete
Schr\"odinger operator. This is
well-known, e.g., \cite{CTT, KL}. Let  $G=(\Er, \Vr, m, \theta)$ be a
weighted graph and  $V:\Vr\to\R$ be a potential. We set the
Schr\"odinger operator
$\Hc:=\Hc|_{\Cc_c(\Vr)}:=(\Delta_{\Er,   \theta}+V(Q))|_{\Cc_c(\Vr)}$. The
domain of its adjoint is given by 
\begin{eqnarray*}
  \Dc(\Hc^*)=\Big\{f\in \ell^2(\Vr, m^2),
  x\mapsto \frac{1}{m^2(x)} \sum_{y\in \Vr} \Er(x,y)
   (f(x)- e^{\rmi \theta_{x,y}}f(y))+V(x)f(x)\in\ell^2(\Vr, m^2)\Big\}. 
\end{eqnarray*}
Then, given $f\in\Dc(\Hc^*)$, one has:
\begin{eqnarray*}
  \left(\Hc^*f\right)(x)=\frac{1}{m^2(x)}\sum_{y\sim 
      x}\Er(x,y)(f(x)-e^{\rmi \theta_{x,y}}f(y))+V(x)f(x),\text{ for all }x\in V.
\end{eqnarray*}
By definition, the operator $\Hc$ is \emph{essentially
self-adjoint} if its closure is equal to its adjoint. Recall that a
symmetric operator is always closable since its adjoint is densely
defined.   

\subsection{Essential self-adjointness}\label{s:esssa}

Before talking about essential self-adjointness, we  deal with
the trivial case, the boundedness of the Laplacian and refer to
\cite{KL} for more discussions in the setting of Dirichlet forms
and $\ell^p$ spaces:

\begin{proposition}\label{p:bounded}
Let $G=(\Vr,  \Er, m, \theta)$ be a weighted graph. One has that
$\Delta_{\Er, \theta}$ is bounded if and only if $d_G(Q)$ is bounded.
\end{proposition} 
\proof 
First, \eqref{e:majo} gives one direction. On the other hand, $\langle
\bone_{\{x\}}, \Delta_{\Er, \theta}\, \bone_{\{x\}}\rangle=d_G(x)$, for all $x\in
\Vr$. \qed 

Essential self-adjointness of the discrete Laplacian was proved in many
situations by J\o rgensen (see \cite{Jor08} and references therein). 
In \cite{Woj} Wojciechowski proves that every discrete Laplacian is
essentially self-adjoint on simple graphs. This result was
independently proved in \cite{Jor08} but the proof was incomplete (see
\cite{JP}).  An alternative proof of this result can be 
found in  \cite{Web} where one uses the maximum principle. Similar
ideas are found in \cite{KL}, where one generalizes this fact to some
weighted graphs by studying Dirichlet forms. Then come the works of   
\cite{Tor, Ma} for weighted graphs which are metrically complete, see
also \cite{CTT} for the non-metrically complete case. Some other
criteria, based on commutators, are given in \cite{GS} (see also
\cite{Gol}). Finally for the magnetic case, we mention the works
\cite{CTT2, Mil, Mil2, MiTr}. We point out that in the older work of \cite{Aom}
ones gives some  characterization of possible self-adjoint extensions
of a weighted discrete Laplacian in the limit point/circle
spirit in the case of trees. More recently, in \cite{GS}, the question
of the deficiency indices is discussed. We point out that, in the
latter, one can consider potentials that tend to $-\infty$.

We now improve a self-adjointness criteria given in \cite{KL} and
extend it in two directions:  we allow magnetic operators and
potentials that are unbounded from below. 

\begin{proposition}\label{p:essaa}
Let $G=(\Vr,  \Er, m, \theta)$ be a weighted graph, $V:\Vr\to \R$ and  
$\gamma>0$. Take $\lambda\in \R$ so that 
\begin{align}\label{e:essaa0}
\{x\in \Vr, \lambda+d_G(x)+ V(x)=0\}=\emptyset.
\end{align}
Suppose that, for any $(x_n)_{n\in
  \N}\in\Vr^\N$, such that the weight $\Er(x_n, x_{n+1})>0$
for all $n\in \N$, the property 
\begin{align}\label{e:essaa1}
\sum_{n\in \N} m^2(x_n) a_n=\infty, \mbox{ where }
a_n:=\prod_{i=0}^{n-1}\left(\left(\frac{\gamma}{d_G(x_i)}\right)^2+ \left(1+
    \frac{\lambda+V(x_i)}{d_G(x_i)}\right)^2\right) 
\end{align}
holds true. Then, the operator $\Hc:=(\Delta_{\Er,    \theta}+V(Q))|_{\Cc_c(\Vr)}$ is
essentially   self-adjoint.    
\end{proposition}
First, note it is always possible to find a $\lambda$ fulfilling
\eqref{e:essaa0}, as $\Vr$ is countable.
Our technique relies  on an improvement of \cite[Theorem 1.3.1]{Woj}. 

\proof  Let $f\in
\Dc(\Hc^*)\setminus \{0\}$ such that $\Hc^*f+ (\gamma \rmi +\lambda) f=0$ or
$\Hc^*f+ (-\gamma\rmi +\lambda) f=0$. We get easily: 
\[|f(x)|\leq \frac{1}{m^2(x)}\sum_{y\in \Vr} 
\frac{\Er(x,y)}{\sqrt{\gamma^2+(\lambda+d_G(x)+V(x))^2} }|f(y)|.\]
We derive:
\begin{align}\label{e:fy}
\max_{y\sim x}|f(y)|^2\geq \left(\left(\frac{\gamma}{d_G(x)}\right)^2+ 
  \left(1+\frac{\lambda+V(x)}{d_G(x)}\right)^2\right)  |f(x)|^2, \mbox{ for all
} x,y\in\Vr, \mbox{ so that } \Er(x,y)\neq 0.
\end{align}
Now, since $f\neq 0$, there is $x_0\in \Vr$ such that $f(x_0)\neq
0$. Therefore, inductively, we obtain a sequence $(x_n)_{n\in \N}\in
\Vr^\N$ such that $\Er(x_n, x_{n+1})>0$, for all $n\in \N$, and so that
\eqref{e:fy} holds for $y=x_{n+1}$ and $x=x_n$. Hence, we get 
\[\sum_{n=0}^N m^2(x_n) |f(x_n)|^2\geq \sum_{n=0}^N m^2(x_n)
\prod_{i=0}^{n-1}\left(
  \left(\frac{\gamma}{d_G(x_i)}\right)^2+\left(1+\frac{\lambda+V(x_i)}
    {d_G(x_i)}\right)^2\right)|f(x_0)|^2.\]    
By letting $N$ go to infinity and remembering \eqref{e:essaa1}, we
obtain a contradiction of the fact that $f\in
\ell^2(\Vr, m^2)$. We conclude with the help of \cite[Theorem X.1]{RS}.\qed

In \cite{KL}, the hypothesis is stronger, i.e., they take $a_n=1$, 
do not consider magnetic fields, and consider potentials that are
bounded from below. We provide two examples which were not covered.

\begin{corollary}\label{c:essaa}
Let $G=(\Vr,  \Er, m, \theta)$ be a weighted
graph, $\varepsilon>0$, $C>0$,  and $V:\Vr\to \R$ such that
\[V(x)\geq -(1-\varepsilon )d_G(x)-C,\]
for all $x\in \Vr$.  
Suppose that, for any $(x_n)_{n\in
  \N}\in\Vr^\N$, such that the weight $\Er(x_n, x_{n+1})>0$
for all $n\in \N$, 
\begin{align}\label{e:cessaa1}
\sum_{n\in \N} m^2(x_n)=\infty 
\end{align}
holds true. Then, the operator  $\Hc:=(\Delta_{\Er,    \theta}+V(Q))|_{\Cc_c(\Vr)}$ is
essentially   self-adjoint. 
\end{corollary} 
Note that the result is optimal as one cannot take $\varepsilon=0$.
Indeed there are simple graphs on which the adjacency matrix $\Ac:=
d_G(Q)- \Delta_{\Er, 0}$ is not essentially self-adjoint on
$\Cc_c(\Vr)$, see \cite{GS} and references 
therein.

\begin{corollary}\label{c:essaa2}
Let $G=(\Vr,  \Er, m, \theta)$ be a weighted
graph, $V:\Vr\to \R$, and $\gamma>0$. 
Suppose that, for any $(x_n)_{n\in
  \N}\in\Vr^\N$, such that the weight $\Er(x_n, x_{n+1})>0$
for all $n\in \N$, 
\begin{align}\label{e:cessaa21}
\sum_{n\in \N} m^2(x_n) \prod_{i=0}^{n-1}\left(
  \frac{\gamma}{d_G(x_i)}\right)=\infty 
\end{align}
holds true. Then, the operator $\Hc:=(\Delta_{\Er,    \theta}+V(Q))|_{\Cc_c(\Vr)}$ is
essentially   self-adjoint. 
\end{corollary}
We stress that in the latter, we make no hypothesis on growth of $V$.

\subsection{Min-max principle}\label{s:min-max}
We recall some well-known results. We refer to \cite{RS}[Chapter XIII.1]
for more details and to \cite{RS}[Chapter XIII.15] for more
applications. We start with the form-version of the standard
variational characterization of the $n$-th eigenvalue.

\begin{theorem}\label{t:min-max}
 Let $A$ be a non-negative self-adjoint operator with form-domain
$\Dc(A^{1/2})$. For all $n\geq 1$, we define:
\begin{align*}
\mu_n(A):=\sup_{\varphi_1, \ldots, \varphi_{n-1}}\inf_{\psi\in [\varphi_1,
  \ldots, \varphi_{n-1}]^\perp} \langle \psi,
A \psi \rangle, 
\end{align*}
where $[\varphi_1,  \ldots, \varphi_{n-1}]^\perp = \{\psi\in \Dc(A^{1/2}),$ so
that $ \|\psi\|=1$ and   $\langle \psi, \varphi_i\rangle=0, $ with $
i=1, \ldots, n-1\}$. Note that 
  $\varphi_i$ are not required to be linearly independent. 

If $\mu_n$ is (strictly) below the essential spectrum of $A$, it is the
$n$-th eigenvalue, counted with multiplicity, and we have:
\begin{align*}
\dim \ran\, \bone_{[0, \mu_n(A)]}(A)=n. 
\end{align*}

Otherwise, $\mu_n$ is the
  infimum of the essential spectrum. Moreover,
  $\mu_{j}=\mu_n$, for all $j\geq n$ and there  are at most $n-1$
  eigenvalues, counted with multiplicity, below the essential
  spectrum. In that case, 
\begin{align*}
\dim \ran\, \bone_{[0,  \mu_n(A)+
  \varepsilon]}(A)= +\infty, \mbox{ for all }\varepsilon>0.
\end{align*}
\end{theorem}
 
\begin{remark}\label{r:multi}
One has a priori no control on the multiplicity of the
(possible) eigenvalue which is at the bottom of the essential spectrum.
\end{remark} 

This ensures the following useful criteria.

\begin{proposition}\label{p:compa}
Let $A,B$ be two non-negative self-adjoint operators. Suppose that
\[\Dc(A^{1/2})\supset \Dc(B^{1/2}) \mbox{ and } 0\leq \langle
\psi, A\psi\rangle 
\leq  
\langle \psi, B \psi\rangle, \]
for all $\psi \in \Dc(B^{1/2})$. Then one has $\inf \sigma_{\rm ess} (A)\leq
\inf \sigma_{\rm ess} (B)$, 
\[\mu_n(A)\leq \mu_n(B), \mbox{ for all }
n\geq 0,\] 
and 
\begin{align}\label{e:N}
N_\lambda(A) \geq N_\lambda(B), 
\mbox{ for } \lambda \in [0,\infty)\setminus \{\inf \sigma_{\rm ess} (B)\},
\end{align}
where $N_\lambda(A):=\dim \ran\, \bone_{[0, \lambda]}(A)$.  

In particular, if $A$ and $B$ have the same
form-domain, then $\sigma_{\rm ess} (A)=\emptyset$ if and only if
$\sigma_{\rm ess} (B)=\emptyset$. 
\end{proposition} 
We refer to \cite{LSW} for some applications of the last statement, in the
context of the absence of the essential spectrum of differential
operators.  Because of Remark \ref{r:multi}, we stress that we have to
remove ``$\inf \sigma_{\rm ess} (B)$'' in  \eqref{e:N}, as the min-max
principle does not decide whether or not the $B$ has an eigenvalue at this
energy. 

\proof  Theorem \ref{t:min-max} permits us to conclude for the
first part. Supposing now they have the same form-domain, by the uniform
boundedness principle, there are $a,b>0$ such that:
\[\langle \psi, A\psi\rangle \leq a \langle \psi, B \psi\rangle +
b\|\psi\|^2 \mbox{ and } \langle \psi, B\psi\rangle \leq a \langle
\psi, A \psi\rangle + b\|\psi\|^2\]
for all $\psi \in \Dc(A^{1/2})= \Dc(B^{1/2})$. By using  the previous
statement twice we get the result.  \qed

We now turn to a criteria of compactness. We recall that a compact
and self-adjoint operator $A$ is in the $p-$Schatten class, for $p\in
[1, \infty)$, if $\tr(|A|^p)<\infty$.
 
\begin{proposition}\label{p:compact}
Let $\Hr$ be a Hilbert space  and $A,B$ be two bounded self-adjoint
operators. Suppose that $B$ is compact and 
\begin{eqnarray*}
|\langle f, A f\rangle|\leq \langle f, B f\rangle, \mbox{ for
  all }f\in \Hr. 
\end{eqnarray*}  
Then $A$ is also compact. Moreover, given $p\in [1, \infty)$, if $B$
is in the $p-$Schatten class, so is $A$. 
\end{proposition} 
\proof First note that $-B\leq A$ and that $-B\leq -A$ in the form
sense. By Proposition \ref{p:compa}, this implies that the essential
spectrum of $A$ is $\{0\}$. Finally, as $A$ is self-adjoint, we infer
that $A$ is compact. Now assume that $B$ is in the $p-$Schatten class.
The min-max principle, applied to $A$ and $-A$, guarantees that
$\tr(|A|^p)\leq 2^p\,\tr(B^p)$. \qed

Using \eqref{e:majo}, we give a straightforward
consequence of Propositions \ref{p:compa} and \ref{p:compact}:  

\begin{corollary}\label{c:essDeltanec}
Let $G=(\Vr,  \Er, m, \theta)$ be a weighted graph. One has:
\[\inf \sigma_{\rm ess}(\Delta_{\Er, \theta})\leq 2 \inf \sigma_{\rm
  ess}(d_G(Q)) \mbox{ and } N_\lambda(\Delta_{\Er, \theta})\geq
N_\lambda(2d_G(Q)),\]
for all $\lambda \in [0, \infty)\setminus \{\inf \sigma_{\rm ess} (d_G(Q))\}$. 
In particular, if $0\in \sigma_{\rm ess} (d_G(Q))$, then $0\in \sigma_{\rm ess}
(\Delta_{\Er, \theta})$ and if $\Delta_{\Er, \theta}$ has compact
resolvent then, $d_G(Q)$ has also compact resolvent. In other words,
one has that $\sigma_{\rm ess} (d_G(Q))\neq 
\emptyset$, then  $\sigma_{\rm ess}  (\Delta_{\Er, \theta})\neq
\emptyset$. 

Moreover, if $d_G(Q)$ is compact, then $\Delta_{\Er, \theta}$ is also
compact. 
\end{corollary} 

We mention, that for a simple graph, $d_G(Q)\geq 1$ and it is
not a compact operator.

\section{Surrounding the Laplacian}
\subsection{A Hardy inequality}\label{s:hardy}
We start with a remark about
bi-partite graphs:

\begin{proposition}\label{p:bipartite} 
Given a bi-partite graph $G=(\Vr,  \Er, m, \theta)$ and a function
$V:\Vr\to [0,\infty)$, the following
assertions are equivalent:
\begin{align}
\label{e:bip1}
\langle f,
(d_G(Q)-V(Q)) f\rangle \leq& \langle f, \Delta_{\Er, \theta} f  \rangle,
&\mbox{ for } f\in 
\Cc_c(\Vr),
\\
\label{e:bip2}
\langle f, \Delta_{\Er, \theta} f  \rangle \leq& \langle f,
(d_G(Q)+V(Q)) f\rangle, &\mbox{ for } f\in
\Cc_c(\Vr),
\\
\label{e:bip3}
|\langle f, \Ac_{\Er, \theta} f\rangle |\leq& \langle f, V f\rangle,
&\mbox{ for } f\in \Cc_c(\Vr),
\end{align}
where $\Ac_{\Er,\theta}$ is the magnetic adjacency matrix defined by
\begin{align}\label{e:adjma}
(\Ac_{\Er,\theta}f)(x):= \frac{1}{m^2(x)}\sum_{y\in \Vr} \Er(x,y) e^{\rmi
  \theta_{x,y}}f(y),
\end{align} 
for  $f\in \Cc_c(\Vr)$ and $x\in \Vr$.
\end{proposition} 
\proof Set $Uf(x):=
(-1)^{|x|}f(x)$. Note that $U^2=\Id$ and $U^{-1}=U^*=U$. Notice now
that on $\Cc_c(\Vr)$
\[U^{-1}\Ac_{\Er,\theta}
U= -\Ac_{\Er,\theta} \mbox{ and } \Delta_{\Er,   \theta}= d_G(Q) -
\Ac_{\Er, \theta}.\] 
We start with \eqref{e:bip1} and rewrite it as follows: 
$\langle f,\Ac_{\Er, \theta} f\rangle\leq \langle f, V(Q) f\rangle$, for
$f\in \Cc_c(\Vr)$.  Applying this to $Uf$, we infer immediately \eqref{e:bip3}.
We start now from \eqref{e:bip3}. We get:
\[\langle f, \Delta_{\Er, \theta}\, f  \rangle = \langle f, (d_G(Q) - \Ac_{\Er,
  \theta}) f  \rangle\geq \langle f, (d_G(Q) - V(Q)) f  \rangle\] 
for  $f\in \Cc_c(\Vr)$. In the same way, we have: \eqref{e:bip2}
is equivalent to \eqref{e:bip3}. 
\qed

We now turn to the key estimate of this paper and prove Proposition
\ref{p:hardy}. First, given $m_0:\Vr\to (0, \infty)$, we mention that
a Laplacian in a certain $\ell^2(\Vr, m^2)$ is unitarily equivalent to
a Schr\"odinger operator in any other  $\ell^2(\Vr, m^2_0)$. This has  
already been noticed before, e.g., \cite{CTT, HK}.

\begin{proposition}\label{p:uni}
Let $G=(\Vr,  \Er, m, \theta)$ and $G_0=(\Vr,  \Er_0, m_0, \theta)$ be
two weighted graphs. Then the Friedrichs extension of $\Delta_{\Er_0, \theta}$,
acting in  $\ell^2(\Vr, m^2_0)$, is unitarily equivalent to  that
of $\Delta_{\Er, \theta}+ V(Q)$, in $\ell^2(\Vr, m^2)$, where
\begin{eqnarray*}
\Er_0(x,y):= \Er(x,y) \frac{m_0(x) m_0(y)}{m(x) m(y)}\quad \mbox{ and } \quad
V(x):=\frac{1}{m_0^2(x)}\sum_{y\in \Vr} \Er(x,y)  \left(1
 - \frac{m(y)}{m(x)}\frac{m_0(x)}{m_0(y)}\right). 
\end{eqnarray*} 
\end{proposition} 
\proof Consider the unitary map $U:\ell^2(G,m^2)\rightarrow
\ell^2(G, m_0^2)$ given by
$(Uf)(x)=(m(x)/m_0(x))f(x)$. Straightforwardly, 
using \eqref{e:Deltam}, one has
\begin{align}\label{e:UU}
U^{-1}\Delta_{\Er_0, \theta}\,U f : =  (\Delta_{\Er, \theta}+ V(Q)) f
\end{align}
for all $f\in \Cc_c(\Vr)$. Then it holds true for the closures of
$\Delta_{\Er_0, \theta}|_{\Cc_c}$ and $(\Delta_{\Er, \theta}+
V(Q))|_{\Cc_c}$. Note now that \eqref{e:UU} also holds for the adjoints of
the last two operators. Therefore, the Friedrichs extensions $\Delta_{\Er_0,
  \theta}$ and $\Delta_{\Er, \theta}+ V(Q)$ are unitarily
equivalent.\qed

The novelty in the Hardy inequality \eqref{e:hardy} is more in the
point of view. Rather than studying $\Delta_{\Er,\theta}$ in $\ell^2(\Vr,
m^2)$ with the help of a simpler $\ell^2(\Vr,m^2_0)$ (where typically
$m_0=1$), we study it with the help of all other weighted
spaces. The applications we consider in this paper are also new. 

\proof[Proof of Proposition \ref{p:hardy}]
On $\ell^2(G,m^2)$, we consider the quadratic form:
\[\Qr_{\Er, \theta}(f,f):= \frac{1}{2}\sum_{x,y\in \Vr} \tilde \Er(x,y)|f(x)-e^{\rmi
  \theta_{x,y}}f(y)|^2\geq 0, \mbox{ for 
} f\in \Cc_c(V),\]
where 
\[\tilde \Er(x,y):=\Er (x,y) \frac{m(x) m(y)}{m_0(x) m_0(y)}.\]
Note that $\Qr(f,f)= \langle f, H_{m^2}
f\rangle_{\ell^2(\Vr,m^2)}$, where 
\[(H_{m^2} f)(x)= \frac{1}{m^2(x)}\sum_{y\sim x} \tilde \Er(x,y)
(f(x)-e^{\rmi \theta_{x,y}}f(y)),  \mbox{ for } f\in \Cc_c(V).\] 
Consider now the unitary map $U:\ell^2(G,m^2)\rightarrow
\ell^2(G, m_0^2)$ given by $(Uf)(x)=(m(x)/m_0(x))f(x)$. Set $H : = U
H_{m^2} U^{-1}$ 
We have:
\begin{align*}
(Hf)(x)&= \frac{1}{m^2(x)}\sum_{y\in \Vr} \tilde \Er(x,y) f(x) -
\frac{1}{m_0^2(x)}\sum_{y\in \Vr} \Er(x,y) e^{\rmi \theta_{x,y}}f(y)
\\
&= \frac{1}{m_0^2(x)}\sum_{y\in \Vr} \Er(x,y) (f(x)- e^{\rmi \theta_{x,y}}f(y))
-V_m(x)f(x).  
\end{align*}
Since $H\geq 0$, we obtain \eqref{e:hardy}. The rest of the statement is
given by Proposition \ref{p:bipartite}. \qed

In this paper, we are interested in minorating the Laplacian with the
help of the weighted degree. Before proving the main result for
simple trees in the next section, we start by giving some negative
answers. First, it is obvious that one cannot find a non-negative
$V_m$ in \eqref{e:hardy} if the graph is finite. This continues to
hold true for some infinite graphs.  

\begin{proposition}\label{p:V=0}
Let $G=(\Vr,  \Er, m, \theta)$ be a weighted graph and take
$V:\Vr\rightarrow \R$ so that we have: $0\leq V(Q)\leq \Delta_{\Er, \theta}$, in the form 
sense on $\Cc_c(\Vr)$. If the constant function $1$ is in
$\Dc((\Delta_{\Er,   \theta})^{1/2})$ then $V=0$. 
\end{proposition} 
If one supposes that $1\in \ell^2(\Vr, m^2)$ and $\Delta_{\Er, \theta}$ is
essentially self-adjoint on $\Cc_c(\Vr)$, then $1$ is in the form domain
$\Dc((\Delta_{\Er,   \theta})^{1/2})$, as $1\in  \Dc((\Delta_{\Er,
  \theta}|_{\Cc_c})^*)$.  

\proof By construction of the Friedrichs extension, there is $f_n\in 
\Cc_c(\Vr)$, so that $f_n$ tends to $1$ in the graph norm of
$\Delta_{\Er,\theta}^{1/2}$. Moreover, 
\[ 0\leq \sum_{x\in \Vr} m^2(x)V(x)|f_n(x)|^2= \langle f_n, V f_n \rangle
\leq \langle f_n, \Delta_{\Er,\theta} f_n\rangle \longrightarrow 0,\]
as $n$  goes to $\infty$. Now since $V^{1/2}(Q)$ is closed, the closed
graph Theorem gives that the l.h.s.\ tends to $\sum_{x\in \Vr}
m^2(x)V(x)$. Therefore $V=0$. \qed 

We now treat the impossibility of minorating by using the
weighted degree.  We say that $\{K_n\}_{n\in \N}$ is a
\emph{filtration} of a graph $(\Vr,  \Er, m, \theta)$ if $K_n\subset
\Vr$ is finite, $K_n\subset K_{n+1}$, for all $n\in \N$, and such
that $\cup_{n\in   \N}K_n=\Vr$. Recall that given a subset $K\in \Vr$, 
$\partial K:=\{x\in K,$ there is $y\sim x$ such that
$y\notin K\}$. 

\begin{proposition}\label{p:negans}
Let  $G=(\Vr,  \Er, m, 0)$ be a weighted graph.

\begin{enumerate}
\item Assume that $\alpha(G)=0$. Then, there is no $a>0$ such that 
\begin{align}\label{e:negans1}
\langle f, d_G(Q) f\rangle\leq a \langle f,
\Delta_{\Er, 0}\, f\rangle, \mbox{ for } f\in\Cc_c(\Vr). 
\end{align} 
\item Suppose that there is a filtration $\{K_n\}_{n\in \N}$ associated with
$G$, such that:
\begin{align}\label{e:negans4}
\limsup_{n\rightarrow +\infty}\frac{\sum_{x\in K_n\setminus \partial K_n}
  m^2(x)d_G(x)}{\sum_{x\in   K_n} m^2(x)}= +\infty,
\end{align}
then there is no $a\in [0,1]$ and $b\geq 0$ such that
\begin{align}\label{e:negans2}
\langle f, d_G(Q) f\rangle\leq a\langle f, \Delta_{\Er,0}\, f
\rangle + b\|f\|^2, \mbox{ for all } f\in \Cc_c(\Vr).
\end{align}
In particular, \eqref{e:negans4} is fulfilled if $1\in \ell^1(\Vr,
m^2)$ and $\sup_{x\in \Vr} m^2(x)d_G(x)=+ \infty$. 
\end{enumerate}  
\end{proposition} 
\proof Take $K\subset \Vr$ finite and use \eqref{e:negans1} with
$f=\bone_{K}$. The statement follows from the definition of $\alpha(G)$. 

For the second
assertion, consider $K=K_n$ and note that 
\eqref{e:negans2} implies that
\[\sum_{x\in
  K_n\setminus \partial K_n} m^2(x)d_G(x) \leq b \sum_{x\in K_n} m^2(x).\]
By letting $n$ go to infinity, we obtain a contradiction with
\eqref{e:negans4}. \qed

\subsection{The case of trees}\label{s:hardytree}
We now turn to a minoration of the magnetic Laplacian and present it in the
context of weighted trees. Perturbation theory will be considered in
Section \ref{s:pertu}. We fix some notation. 

We first recall that a \emph{tree} is a connected graph $G=(\Er,\Vr,
m, \theta)$ such that for each edge $e\in \Vr\times \Vr$ with
$\Er(e)\neq 0$ the graph $(\tilde  \Er, \Vr, m, \theta|_{\tilde
  \Er})$, with $\tilde 
\Er := \Er \times 1_{\{e\}^c}$, i.e., with the edge $e$ 
removed, is disconnected. It is convenient to choose a root in the
tree. Due to its structure, one can take any point of $\Vr$. As in the
definition of $|x|$, we choose $\omega$ to be the root. We define
the \emph{sphere} $S_n$ by 
\[S_{-1}=\emptyset, S_0:=\{\omega\}, \mbox{ and }
S_{n}:=\{x\in \Vr, |x|=n\}.\]
%
  Given $n\in\N$, $x\in S_n$, and
$y\in\Nr_G(x)$,  one sees that $y\in S_{n-1}\cup S_{n+1}$. 
We write $x\rightsquigarrow y$ and say that $x$ is a \emph{son}
of~$y$, if $y\in S_{n-1}$, 
while we write $x{\reflectbox{$\rightsquigarrow$}}  y$ and say that $x$ is
a \emph{father} of~$y$, 
if $y\in S_{n+1}$.
Notice that $\omega$ has no father.
Given $x\ne\omega$, note that there is a unique $y\in V$ with $x\rightsquigarrow y$,
i.e., everyone apart from $\omega$ has one and only one father.
We denote the father of $x$ by $\father x$.

We turn to the proof of Theorem \ref{t:tree}. We strengthen slightly
the result by working under the
hypothesis 
\begin{align}\label{e:hypC0}
C_0:=\sup_{x\in \Vr}
  d^{\varepsilon_0 -1}_{G_\circ}(x) \Er_\circ(x,\father x)m^{-2}(x)
  <\infty \mbox{ and } C_1:=\inf_{x\in \Vr} d_{G_\circ}(x)>0.
\end{align}
instead of \eqref{e:hypC0old}. 

\begin{remark}\label{r:hypC0}
We mention that because of  part (b) of Proposition
\ref{p:negans}, one cannot solely suppose  that $G$ is a weighted
tree. This is why a relation between $m$ and $\Ec$ has to be
assumed. Whereas the condition on $C_1$ is easy to check, the
condition on $C_0$ seems more technical. In order to understand it
better, we point out that it  would be automatically fulfilled if we
could take $\varepsilon_0=0$. The optimality of \eqref{e:hypC0}
remains open. 
\end{remark} 

\proof[Proof of Theorem \ref{t:tree}] Before going into perturbation
theory, we focus on the left hand side of \eqref{e:tree1} for
$G_\circ$ and $V=0$, i.e., for $\Delta_{\Er_\circ,
  \theta}$ instead of $\Hc_\Fr$.  Take $\eta>0$. We define:
\[\tilde m(\omega):=1 \mbox{ and } \tilde m(x):=
\eta\,  \tilde m(\father  x)\frac{m(x)}{m(\father
  x)}d_{G_\circ}^{-\varepsilon_0/2}(x), 
\mbox{ for all } x\in  
\Vr\setminus\{\omega\}.\]  
With this definition and $V_{\tilde m}$ as in \eqref{e:Vm}, we obtain:
\begin{align}\nonumber
\frac{V_{\tilde m}(x)}{d_{G_\circ}(x)}&=
1-\frac{1}{d_{G_\circ}(x)m^2(x)}\left(\Er_\circ(\father x, x) \frac{\tilde m(\father
    x)}{\tilde m(x)}\frac{m(x)}{m(\father x)} +  
  \sum_{y\rightsquigarrow  x} \Er_\circ(y,   x) \frac{\tilde m(y)}{\tilde
    m(x)}\frac{m(x)}{m(y)} \right)  
\\ \nonumber
&= 1 -  \frac{1}{\eta}\frac{1}{d_{G_\circ}^{(1-\varepsilon_0/2)}(x)m^2(x)}
\Er_\circ(\father x,x)
- \frac{\eta}{d_{G_\circ}(x)m^2(x)}\sum_{y\rightsquigarrow  x}
\Er_\circ(y,x) d_{G_\circ}^{-\varepsilon_0/2 }(y) \\\label{e:key2}
& \geq 1 - \eta C_1^{-\varepsilon_0/2 } -
\frac{1}{\eta}d_{G_\circ}^{-\varepsilon_0/2 }(x) C_0, 
\end{align}
for all $x\in \Vr$ and where $C_0$ and $C_1$ are given by
\eqref{e:hypC0}. Now note that for all
$\tilde \eta>0$ there is 
$c_{\tilde \eta}>0$ such that $t^{1-\varepsilon_0/2}\leq \tilde \eta t + c_{\tilde \eta}$,
for all $t\in [0,\infty)$. Therefore by Proposition 
\ref{p:hardy} applied to $G_\circ$, for all $\varepsilon>0 $, there
is $c_\varepsilon >0$ such that: 
\begin{align}
\nonumber
\langle f, \Delta_{\Er_\circ, \theta} f\rangle & \geq
 \langle d_{G_\circ}^{1/2}(Q) f,\left(1- \eta C_1^{-\varepsilon_0/2 }-  
C_0/d_{G_\circ}^{\varepsilon_0 /2}(Q)\eta\right)  d_{G_\circ}^{1/2}(Q)
f\rangle 
\\ 
\label{e:tree1'}
&\geq (1-\varepsilon ) \langle f, d_{G_\circ}(Q) f\rangle -
c_\varepsilon \|f\|^2,
\end{align}
for all $f\in \Cc_c(\Vr)$.
This gives \eqref{e:tree1} for $G_\circ$ and $V=0$, where the second
inequality is obtained by applying Proposition \ref{p:bipartite}.  

Next, the equality of the
domains of  the forms follows immediately and so the
essential spectrum of $\Delta_{\Er_\circ, \theta}$ is empty if and only if 
$\lim_{|x|\rightarrow   \infty} d_{G_\circ}(x)= +\infty$. Finally, we use
twice Proposition \ref{p:compa} with the double inequality
\eqref{e:tree1}. This yields: 
\[ 1- \varepsilon \leq  \liminf_{N \to
  \infty}\frac{\lambda_N(\Delta_{\Er,  \theta})}{\lambda_N(d_{G_\circ}(Q))}
\leq \limsup_{N \to  \infty}\frac{\lambda_N(\Delta_{\Er,
 \theta})}{\lambda_N(d_{G_\circ}(Q))}\leq 1+ \varepsilon.\] 
By letting $\varepsilon$ go to zero we obtain the  asymptotic
\eqref{e:asymp} for $\Delta_{\Er_\circ ,  \theta}$. 

We finish with $\Hc_\Fr=\Delta_{\Er,  \theta}+V(Q)$ by perturbing
$\Delta_{\Er_\circ, \theta}$. Thanks to  
\eqref{e:tree0}, Propositions \ref{p:pertutree} and \ref{p:pertutree2}
end the proof. \qed

\begin{remark}
 It is important to note that $\varepsilon$ has to be positive in
 \eqref{e:tree1}. For instance, on a simple tree, if $\varepsilon$
 was equal to $0$,  this  would imply that the magnetic adjacency
 matrix was bounded (see 
 Proposition \ref{p:bipartite}). Finally note that, by considering
 delta functions as test functions, the magnetic  adjacency matrix is
 bounded if and only if the weighted degree is bounded.  
\end{remark}

\begin{remark}
We could not provide an example of a Schr\"odinger operator for which we
can compute the asymptotic of eigenvalues and which is not essentially
self-adjoint on $\Cc_c(\Vr)$.  
\end{remark} 

\section{Comparison of domains}
\subsection{From form-domain to domain}\label{s:dom}
Once one knows that the form-domains are equal, the next question is
to guarantee that the domains are equal. We present an approach
by commutators. Some subtleties arise since we have to deal with the
square root of the Laplacian. 

We first prove the invariance of the domain of $\Delta_{\Er, \theta}$
under the $C_0-$group $\{e^{\rmi t d_G(Q)}\}_{t\in\R}$ and that
$\Delta_{\Er, \theta}\in \Cc^1(d_G(Q)) $, i.e., the map $t\mapsto e^{-\rmi t
  d_G(Q)}(\Delta_{\Er, \theta}+1)^{-1}e^{\rmi t
  d_G(Q)} $ is strongly $\Cc^1$. We refer to Appendix \ref{a:regu} for
a discussion about the $\Cc^1$ regularity.

\begin{proposition}\label{p:invdom}
Let $G=(\Vr,  \Er, m, \theta)$ be a weighted graph, such that
$\Delta_{\Er, \theta}$ is essentially self-adjoint on $\Cc_c(\Vr)$, 
$\Dc(\Delta_{\Er, \theta}^{1/2})= \Dc(d_G^{1/2}(Q))$, and 
\begin{align}\label{e:invdom}
\sup_{x\in \Vr} \left(\frac{1}{m^2(x)}\sum_{y\in \Vr}
\frac{\Er(x, y)}{\min (\langle  d_G(x)\rangle^{1/2}, \langle
  d_G(y)\rangle^{1/2})} \, |d_G(x)-d_G(y)|\right) <\infty.
\end{align}
Then, one has the invariance
$e^{\rmi t d_G(Q)}\Dc(\Delta_{\Er, \theta}) \subset \Dc(\Delta_{\Er,
  \theta})$, for all $t\in \R$. Moreover, we have that $\Delta_{\Er, \theta}\in
\Cc^1(d_G(Q))$ and that $[\Delta_{\Er, \theta}, d_G(Q)]|_{\Cc_c(\Vr)}$ 
extends to $[\Delta_{\Er, \theta},
d_G(Q)]_\circ\in\Bc\left(\Dc(\Delta^{1/2}_{\Er,   \theta}),
  \ell^2(\Vr, m^2)\right)$. 
\end{proposition} 
\proof We denote by $[\Delta_{\Er, \theta}, e^{\rmi t d_G(Q)}]_\circ$ the
closure of the commutator $[\Delta_{\Er, \theta}, e^{\rmi t
  d_G(Q)}]|_{\Cc_c(\Vr)}$ and by $c_0$ the constant appearing in
\eqref{e:invdom}. We have:
\begin{align}
\nonumber
|\langle f, [\Delta_{\Er, \theta}, e^{\rmi t d_G(Q)}]_\circ\, g  \rangle|&
=   \sum_{x\in \Vr} \sum_{y\in \Vr} \overline{f(x)}\Er(x, y) e^{\rmi \theta_{x,y}}
\left(e^{\rmi   t d_G(x)}-e^{\rmi t d_G(y)}\right) g(y) 
\\ 
\nonumber
&\leq \frac{|t|}{2} \sum_{x\in \Vr} \sum_{y\in \Vr} \frac{\Er(x,
  y)}{\langle d_G(y)\rangle^{1/2}} \,
|d_G(x)-d_G(y)| \left(| f(x)|^2+ |\langle
  d_G(y)\rangle ^{1/2} g(y)|^2\right)
\\
\label{e:invdom1}
&\leq c_0 |t| \left(\left\|
f\right\|^2 + \left\| \langle 
d_G(Q)\rangle ^{1/2}g\right\|^2 \right),
\end{align}
for all $f, g\in \Cc_c(\Vr)$. Therefore, $[\Delta_{\Er, \theta},
e^{\rmi t d_G(Q)}]_\circ$ is bounded from $\Dc(\Delta^{1/2}_{\Er,
  \theta})$ to $\ell^2(\Vr, m^2)$.  On the other hand, 
we get:
\begin{align*}
\|\Delta_{\Er, \theta}\, e^{\rmi t d_G(Q)} f \|&\leq \|e^{\rmi t
  d_G(Q)}\Delta_{\Er, \theta}\, f\| + \|[\Delta_{\Er, \theta}, e^{\rmi t
  d_G(Q)} ] \langle d_G(Q)\rangle^{-1/2} \langle d_G(Q)\rangle^{1/2}f\|
\\
&\leq C \left(\|f\|+ \|\Delta_{\Er, \theta}\, f\| \right),
\end{align*}
for all $f\in \Cc_c(\Vr)$. By density we obtain the invariance of the domain. 

Next, thanks to \eqref{e:invdom1}, we obtain that
\[ \liminf_{t\rightarrow 0^+}\left\| [\Delta_{\Er, \theta}, e^{\rmi t
  d_G(Q)}/t] \right\|_{\Bc(\Dc(\Delta_{\Er, \theta}), \ell^2(\Vr, m^2))} <\infty.\]
Therefore, Theorem \ref{th:abg2} yields that $\Delta_{\Er,
  \theta}\in \Cc^1(d_G(Q))$. Finally, by estimating as in
\eqref{e:invdom1}, we obtain that $[\Delta_{\Er, \theta},
d_G(Q)]_\circ$, belongs to $\Bc\left(\Dc(\Delta^{1/2}_{\Er,   \theta}),
  \ell^2(\Vr, m^2)\right)$. \qed

We turn to the central result of this section.

\begin{theorem}\label{t:equadom}
Let $G=(\Vr,  \Er, m, \theta)$ be a weighted graph, such that
$\Dc(\Delta_{\Er,   \theta}^{1/2})= \Dc(d_G^{1/2}(Q))$, $\Delta_{\Er,
  \theta}$ is essentially self-adjoint on $\Cc_c(\Vr)$, and there is
$\varepsilon>0$ such that
\begin{align}\label{e:equadom}
\sup_{x\in \Vr} \left(\frac{1}{m^2(x)}\sum_{y\in \Vr}
\frac{\Er(x, y)}{\min (\langle  d_G(x)\rangle^{1/2- \varepsilon}, \langle
  d_G(y)\rangle^{1/2-\varepsilon})} \, |d_G(x)-d_G(y)|\right) <\infty.
\end{align}
Take also a potential $V:\Vr\rightarrow \R$, such that:
\begin{align}\label{e:equadom2}
\lim_{|x|\rightarrow +\infty} \frac{V(x)}{d_G(x)+1}=0.
\end{align}
Then, the operator $\Hc:=(\Delta_{\Er,   \theta}+V(Q))|_{\Cc_c(\Vr)}$
is bounded from below. We denote by $\Hc_\Fr$ its Friedrichs
extension. We have: $\Dc(\Hc_{\Fr})= \Dc(d_G(Q))$ and $\sigma_{\rm
ess}(\Hc_{Fr})= \sigma_{\rm
ess}(\Delta_{\Er,   \theta})$.
\end{theorem} 

\begin{remark}
We point out that the hypothesis on the essential self-adjointness is
also necessary, as $d_G(Q)$ is essentially self-adjoint on
$\Cc_c(\Vr)$.  
\end{remark}

\begin{remark}
By taking $\varepsilon=0$ in \eqref{e:equadom}, i.e., under the same
hypothesis as Proposition \ref{p:invdom}, the proof ensures only that
$\Dc(d_G(Q))\subset \Dc(\Delta_{\Er, \theta})$. 
\end{remark} 

\proof The stability of the essential spectrum is ensured by
Proposition \ref{p:pertutree2}. We focus on the domain.
 By Kato-Rellich's Theorem, e.g., \cite[Theorem X.12]{RS}, it is
enough to consider $V=0$.  We assume for the moment that 
\[M:= (\Delta_{\Er,   \theta}+1)^{1/2} [d_G(Q),  (\Delta_{\Er,
  \theta}+1)^{1/2}]_\circ +  [\Delta_{\Er,   \theta},
(d_G(Q)+1)^{1/2}]_\circ (d_G(Q)+1)^{1/2}.\] 
is a bounded operator from $\Dc(d_G^{1-\varepsilon }(Q))$ to
$\Dc(d_G^{1-\varepsilon }(Q))^*$. The $1$ is here to make the square root smooth
over the spectrum. Since the form-domain  of
$\Delta_{\Er,   \theta}$ and $d_G(Q)$ are equal, by the uniform boundedness
principle, we have: there are $a,b>0$
so that \eqref{e:negans2} holds true. By using twice \eqref{e:negans2}
and working in the form  sense on $\Cc_c(\Vr)$, we infer:  

\begin{align}
\nonumber
(\Delta_{\Er,   \theta}+1)^2 &\leq  a (\Delta_{\Er,   \theta}+1)^{1/2} (d_G(Q)+1)
(\Delta_{\Er,   \theta}+1)^{1/2} + (b+1) (\Delta_{\Er,   \theta}+1)
\\
\nonumber
& = a\, (d_G(Q)+1)^{1/2} \Delta_{\Er,  \theta}\, (d_G(Q)+1)^{1/2} + aM
+(b+1) (\Delta_{\Er,   \theta}+1)  
\\
\label{e:equadom3}
& \leq a' d_G^2(Q) + b',
\end{align}
for some $a', b'>0$. By using the fact that $\varepsilon>0$, the
reverse inequality holds true for the same reasons for some $a'\in
(0,1)$ and  $b'>0$. Therefore, the domains are equal. It remains to
prove the  
boundedness of $M$. We start with the r.h.s.\ term. Let $c_1$ be the
constant in \eqref{e:equadom}. We estimate as above:
\begin{align}
\nonumber
|\langle f, [\Delta_{\Er, \theta}, (d_G(Q)+1)^{1/2}]_\circ  \langle 
d_G(Q)\rangle^{1/2} g  \rangle|& 
\\
\nonumber
&\hspace{-2cm}=   \sum_{x\in \Vr} \sum_{y\in \Vr}
\overline{f(x)}\Er(x, y)e^{\rmi \theta_{x,y}}  
((d_G(x)+1)^{1/2}-(d_G(y)+1)^{1/2})
\langle d_G(y)\rangle^{1/2}g(y) 
\\
\nonumber
&\hspace{-2cm} \leq \frac{1}{2} \sum_{x\in \Vr} \sum_{y\in \Vr} \frac{\Er(x,
  y)}{\langle d_G(y)\rangle^{1/2-\varepsilon}} \,
|d_G(x)-d_G(y)| \left(| f(x)|^2+ |\langle
  d_G(y)\rangle^{1-\varepsilon} g(y)|^2\right)
\\
\nonumber
&\hspace{-2cm}
\leq c_1 \left(\left\|
f\right\|^2 + \left\| \langle 
d_G(Q)\rangle ^{1-\varepsilon}g\right\|^2 \right),
\end{align}
for all $f, g\in \Cc_c(\Vr)$. In particular, $[\Delta_{\Er, \theta},
(d_G(Q)+1)^{1/2}]_\circ  \langle  d_G(Q)\rangle^{1/2}$ is bounded from
$\Dc(d_G^{1-\varepsilon}(Q))$ to $\ell^2(\Vr, m^2)$. 

We turn to the second part of $M$ and use Appendix
\ref{a:comm}. Let $\varphi$ be in  $\Sc^{1/2}$ such that 
$\varphi(x)=\sqrt{x}$, for all $x\geq 1$. Since $\Delta_{\Er,
  \theta}$ is non-negative, $\varphi(\Delta_{\Er,
  \theta}+1)= (\Delta_{\Er,\theta}+1)^{1/2}$.  
We cannot use the \eqref{eq:int} directly with $\varphi$ as the
integral does not seem to exist. We proceed as in
\cite{GoleniaJecko}. Take $\cchi_1\in \Cc^\infty_c(\R;\R)$ with values
in $[0,1]$ and being $1$ on $[-1,1]$. Set
$\cchi_R:=\cchi(\cdot/R)$. As $R$ goes to infinity, $\cchi_R$
converges pointwise to $1$. Moreover, $\{\cchi_R\}_{R\in [1,\infty)}$
is bounded in $\Sc^0$. We infer $\varphi_R:=\varphi \cchi_R$ tends
pointwise  to $\varphi$ and that $\{\varphi_R\}_{R\in [1,\infty)}$ is
bounded in $\Sc^{1/2}$. Now, recalling $\Delta_{\Er, \theta}\in
\Cc^1(d_G(Q))$ and \eqref{e:facto}, we obtain 
\begin{align}\nonumber  
\hspace*{-1cm} 
[\varphi_R(\Delta_{\Er, \theta}),
d_G(Q)]_\circ\langle \Delta_{\Er, \theta} \rangle^{-1/2} =& 
\\ \label{e:ligne0}
&\hspace{-4cm}=  \frac{\rmi}{2\pi}\int_\C\frac{\partial\varphi^\C_R }{\partial 
\overline{z}}(z-\Delta_{\Er,
\theta})^{-1} [d_G(Q), \Delta_{\Er, \theta}]_\circ (z-\Delta_{\Er,
\theta})^{-1}\langle \Delta_{\Er, \theta} \rangle^{-1/2} dz\wedge d\overline{z}.   
\\
\nonumber
&  
\hspace{-4cm}=  \frac{\rmi}{2\pi}\int_\C\frac{\partial\varphi^\C_R }{\partial 
\overline{z}}(z-\Delta_{\Er, \theta})^{-1} 
[d_G(Q), \Delta_{\Er, \theta}]_\circ \,  \langle
\Delta_{\Er,\theta}\rangle^{-1/2} (z-\Delta_{\Er,
  \theta})^{-1}
dz\wedge d\overline{z}. 
\end{align} 
By Proposition \ref{p:invdom}, $  [d_G(Q),
\Delta_{\Er,\theta}]_\circ\langle \Delta_{\Er,\theta}\rangle^{-1/2}$
is bounded. Moreover, using \eqref{eq:dg1}  
with $l=2$, we  bound the integrand,
uniformly in $R$, by
$C \langle x\rangle^{-1+1/2 -2}|y|^2 |y|^{-1} |y|^{-1}$, 
for some constant $C$. It is integrable on the domain given by
\eqref{eq:dg3}. By Lebesgue domination, the r.h.s.\ of
\eqref{e:ligne0} has a limit in norm. Note now that the l.h.s., as
form on $\Cc_c(\Vr)$, tends to the operator $[(\Delta_{\Er,
  \theta}+1)^{1/2}, d_G(Q)]\langle \Delta_{\Er,
  \theta}\rangle^{-1/2}$. This gives that
\[-C (d_G(Q)+1)\leq \langle \Delta_{\Er,
  \theta}\rangle^{1/2}[(\Delta_{\Er,
  \theta}+1)^{1/2}, d_G(Q)]_\circ \leq C (d_G(Q)+1)\]
in the form sense on $\Cc_c(\Vr)$ for some constant $C$. 
This ensures the announced boundedness of $M$. \qed

\begin{remark}
Assuming also that for all $\varepsilon>0$, there is $c_\varepsilon \geq
0$ so that \eqref{e:tree1} holds true, one observes that $a'$ can
be arbitrary close to $1$ in \eqref{e:equadom3}. Therefore, using
again the Kato-Rellich's Theorem, one can weaken
\eqref{e:equadom2} in Theorem \ref{t:equadom} and replace it by 
\begin{align}
\limsup_{|x|\rightarrow +\infty} \frac{|V(x)|}{d_G(x)+1}<1.
\end{align}
to obtain $\Dc(\Hc_{\Fr})= \Dc(d_G(Q))$. 
\end{remark} 


Finally, we examine the question of the equality of the domains in the
context of simple trees. We provide  a positive result. 

\begin{proposition}\label{p:treenotdom}
There is a simple tree $T$, such that $\sigma(\Delta_T)=\sigma_{\rm
  ac}(\Delta_T)=[0, \infty)$ and such that $\Dc(\Delta_T)=\Dc(d_T(Q))$. In
particular, $\alpha(T)=0$.  
\end{proposition}
Here ac stands for absolutely continuous. 
\proof 
We start by fixing some notation. Given an \emph{offspring sequence}
$(b_n)_{n\in\N}$, with $b_n\in \N^*$, we associate a simple tree 
with root $\varepsilon$ such that, for all $x\in S_n$, $b_n=\sharp\{y, \father y=x\}$.
\begin{align*}
  \xymatrix{%
    &&\omega\ar@{-}[d]\ar@{-}[dr]&&&& S_0
    \\
    &&{\bullet}\ar@{-}[d]\ar@{-}[dl]\ar@{-}[dll]&
      {\bullet}\ar@{-}[rrd]\ar@{-}[d]\ar@{-}[dr]&&& S_1
    \\
   {\bullet}\ar@{.}[d]&
   {\bullet} \ar@{.}[d]&
   {\bullet}\ar@{.}[d]&
   {\bullet}\ar@{.}[d]&
  {\bullet} \ar@{.}[d]&
  {\bullet} \ar@{.}[d]& S_2
    \\
    &&&&&&
  }
  \\
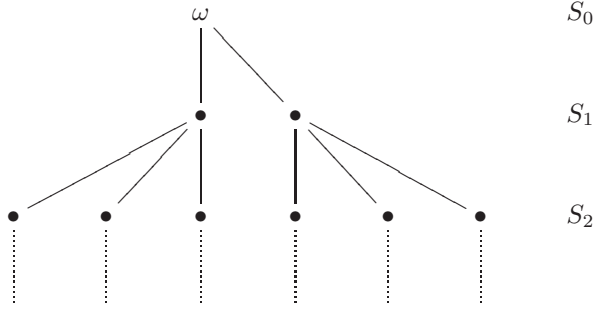
  \text{\it Example of a tree with } b_0=2\, \text{\it and } b_1=3.
  \quad\quad
\end{align*}
We turn to our example and construct some trees $T_n= (\Er_n,
\Vr_n)$.  For $n=1$,  we take $\Vr_1:=\N$, with $\omega_1:=0$ and
with $\Er_1(p,q)=1$ if and only if $|p-q|=1$. For $n\geq 2$, we take
trees that are  $n$-ary after the first generation. For each $n\in 
\N\setminus\{0,1\}$, let $\omega_n$ be the root and set that the
offspring $b(n)_k:=n$, for all $k\in \N\setminus\{0\}$ and 
$b(n)_0:=n-1$. Now take  $T:=(\Er, \Vr)$, where $\Vr:=\cup_{n\in
  \N\setminus\{0, 1\}} \Vr_n$ and $\Er(x,y):=\Er_n(x,y)$, if $x, y\in
\Vr_n$, $\Er(\omega_n,\omega_{n+1}):=1$, and $\Er(x,y):=0$ otherwise. 
\begin{align*}
\begin{array}{ccccccc}
\xymatrix@C=1pc{
{\bullet}\ar@{-}[d]
&
&{\bullet}\ar@{-}[dr] &
&{\bullet}\ar@{-}[dl]
&
{\bullet}\ar@{-}[dr]&{\bullet}\ar@{-}[d]&{\bullet}\ar@{-}[dl]
&
&{\bullet}\ar@{-}[dr]&{\bullet}\ar@{-}[d]
&{\bullet}\ar@{-}[dl]
\\
{\bullet}\ar@{-}[d]
&
&&{\bullet}\ar@{-}[d]
&
&
&
{\bullet}\ar@{-}[drr]&&&&{\bullet}\ar@{-}[dll] &
\\
\omega_1\ar@{--}[rrr]
&
&&\omega_2\ar@{--}[rrrrr]&&
&&
&\omega_3\ar@{--}[rrr]&&&&
\\
} 
\\
\hspace*{-2cm} T_{1} \hspace*{2cm} T_{2} \hspace*{4.1cm} T_3
\\
\\
\mbox{ Graph of } T
\end{array}
\end{align*}
Note that $x\mapsto\sum_{y\sim x} |d_T(x)-d_T(y)|$ has support
contained in $\cup_n\{\omega_n\}$ and takes values in
$\{0,1\}$. Hence, \eqref{e:equadom} is fulfilled. Since the graph
is simple, $\Delta_T$ is essentially self-adjoint on $\Cc_c(\Vr)$ (see
Section \ref{s:esssa}). Therefore, we derive from Theorem
\ref{t:equadom} that the $\Dc(\Delta_T)=\Dc(d_T(Q))$. 

We turn to the spectrum. First, $\sigma(\Delta_{G_1})=\sigma_{\rm
  ac}(\Delta_{G_1})=[0,4]$, where $G_1:=T_1$. This is easy to see by
discrete Fourier transformation, e.g., \cite{AF}. For each $n\geq 2$,
$T$ contains a subtree $G_{n}$ which is $n-$ary and which is connected
to the rest of $T$ by only one edge. It is well-known that
$\sigma_{\rm   ac}(\Delta_{G_{n}}) = [n+1-2\sqrt{n}, n+1+2\sqrt{n}]$, e.g.,
\cite{AF}.  Now, we denote by $\tilde G_i:=T\setminus G_i$, the graph
obtained by removing the only edge that is connecting $G_i$ to the
rest of the graph. Note that $\Delta_{\tilde G_i}-\Delta_{G_i}$ is a
rank one operator, for all $i\geq 1$. Therefore, $\sigma_{\rm
  ac}(\Delta_{G_i})\subset \sigma_{\rm   ac}(\Delta_T)$ for all
$i\geq 1$. Hence, $ \sigma(\Delta_T)=\sigma_{\rm   ac}(\Delta_T)=[0,
\infty)$. Finally note that if
$\alpha(T)>0$, \eqref{e:iso} ensures that $\inf
\sigma(\Delta_T)>0$. This is a contradiction. \qed

\begin{remark}
This construction also provides an example of a graph on which the
adjacency matrix $\Ac_{\Er, 0}$ (see \eqref{e:adjma}) is essentially
self-adjoint on $\Cc_c(\Vr)$ (see \cite[Lemma 2.1]{GS}) and has absolutely
continuous spectrum equal to $\R$. 
\end{remark} 

Finally, we give a negative example.

\begin{proposition}\label{p:treedom}
There is a simple tree $T$, such that $\Dc(\Delta_T)\neq \Dc(d_T(Q))$
and such that the form-domains $\Dc(\Delta_T^{1/2})= \Dc(d_T^{1/2}(Q))$.
\end{proposition}
\proof The second point follows from Theorem \ref{t:tree}. We start by
constructing the star graph $S_n$. Let 
$S_{n}:=(\Er_n, \Vr_n)$ be defined as follows: $\Vr_n:= \{1, \ldots,
n+1\}$ and so that  $\Er_n((1, j)):=1, \forall
j \in \{2, \ldots, n+1\}$  and $\Er_n((j,k)):=0, \forall
j,k \in \{2, \ldots, n+1\}$. Consider now $f_n(x):=1$ on $\Vr_n$. One has:
\begin{align}\label{e:fn2}
\|\Delta_{S_n} f_n\|^2 =0, \quad \|f_n\|^2=n+1, \quad
 \mbox{and} \quad  \|d_{S_n}(Q) f_n\|^2= n(n+1).
\end{align}
Take  $T:=(\Er, \Vr)$, where $\Vr:=\cup_{n\in \N\setminus\{0\}} \Vr_n$ and
$\Er(x,y):=\Er_n(x,y)$, if $x, y\in \Vr_n$, $\Er(x,y):=1$, if
$x=1\in \Vr_n$ and $y=1\in \Vr_{n+1}$, for all $n\geq 1$, and
$\Er(x,y):=0$ otherwise.
\begin{align*}
\begin{array}{ccccccc}
\xymatrix{
{\bullet}\ar@{-}[d]&{\bullet}\ar@{-}[d] &{\bullet}\ar@{-}[dl]&{\bullet}\ar@{-}[d] &{\bullet}\ar@{-}[dl]
&{\bullet}\ar@{-}[dr] &{\bullet}\ar@{-}[d]&{\bullet}\ar@{-}[dl]
\\
{\bullet}\ar@{--}[r]&{\bullet} \ar@{--}[r]&\ar@{--}[r]&{\bullet} \ar@{--}[r]&
\ar@{--}[r]&\ar@{--}[r]&{\bullet} \ar@{--}[r]&
\\
{\bullet}\ar@{-}[u]&{\bullet}\ar@{-}[u]&&{\bullet}\ar@{-}[u]
&{\bullet}\ar@{-}[ul] 
&&{\bullet}\ar@{-}[u] &{\bullet}\ar@{-}[ul]
} 
\\
\hspace*{-1cm} S_{2} \hspace*{1.2cm} S_{3} \hspace*{2.5cm}
S_4\hspace*{3cm} S_5
\\
\\
\mbox{ Graph of } T
\end{array}
\end{align*}
Now if $\Dc(d_T(Q))\supset \Dc(\Delta_T)$, by the uniform boundedness
principle, there are 
constants $a,b>0$, so that
\[ \|d_T(Q) f\|^2 \leq a \|\Delta_T f\|^2 + b\|f\|^2,
\quad \mbox{ for all } f\in \Cc_c(T).\]
This leads to a contradiction with \eqref{e:fn2}, as $\|\Delta_T
f_n\|^2=4$ and $\|d_T(Q) f_n\|^2=(n+2)^2+n$ for $n\geq 3$.\qed

\subsection{The form-domain for  bi-partite graphs}\label{s:dombip} 
As \eqref{e:tree1} holds
true for trees, it is natural to ask the question for bi-partite
graphs.  The answer is no. We start by relating the form-domain of the
magnetic Laplacian with the inequality \eqref{e:reci0}.

\begin{proposition}\label{p:reci}
Let $G=(\Vr,  \Er, m, \theta)$ be a weighted bi-partite graph. Then there is
$a\in (0, 1]$ and $C_a>0$ so that:
\begin{align}\label{e:reci0}
(1-a)\langle f, d_G(Q) f\rangle - C_a \|f\|^2\leq \langle f, \Delta_{\Er,
  \theta} f\rangle \leq    (1+a)\langle f, d_G(Q) f\rangle + C_a \|f\|^2,
\end{align}
for all $f\in \Cc_c(\Vr)$. 
Moreover, one can take some $a<1$ in \eqref{e:reci0} if and only if 
$\Dc((\Delta_{\Er, \theta})^{1/2})= \Dc(d_G^{1/2}(Q))$.
 
Suppose also that 
$\Delta_{\Er, \theta}$ has compact resolvent. Then, 
$d_G(Q)$ has also compact resolvent and, with the same $a$ as in
\eqref{e:reci0}, one has:
\begin{align}\label{e:reci}
 1- a \leq  \liminf_{\lambda \to
  \infty}\frac{N_\lambda(\Delta_{\Er,  \theta})}{N_\lambda(d_G(Q))}
\leq \limsup_{\lambda \to  \infty}\frac{N_\lambda(\Delta_{\Er,
 \theta})}{N_\lambda(d_G(Q))}\leq 1+ a.
\end{align}
\end{proposition} 
We recall that \eqref{e:majo} ensures that $\Dc(d_G^{1/2}(Q))\subset
\Dc((\Delta_{\Er, \theta})^{1/2})$.  

\proof First, note that \eqref{e:reci0} follows from \eqref{e:majo} and 
Proposition \ref{p:bipartite}.  
Suppose that \eqref{e:reci0} holds true for some $a<1$, this gives
immediately that 
$\Dc((\Delta_{\Er,  \theta})^{1/2})=\Dc(d_G^{1/2}(Q))$. Reciprocally,
suppose now  that
$\Dc((\Delta_{\Er,   \theta})^{1/2})= \Dc(d_G^{1/2}(Q))$. 
By the uniform boundedness
principle, there is $a_0, b_0>0$ so that
\[a_0\langle f, d_G(Q) f\rangle - b_0 \|f\|^2\leq \langle f, \Delta_{\Er,
  \theta} f\rangle,\]
for all $f\in \Cc_c(\Vr)$. Using again Proposition \ref{p:bipartite},
\eqref{e:reci0} holds true with $a=1-a_0<1$. 

We turn to the second part and work under the hypothesis that 
$\Delta_{\Er, \theta}$ has compact resolvent. Corollary
\ref{c:essDeltanec} ensures that $d_G(Q)$ has also compact
resolvent. We conclude by using Proposition \ref{p:compa} twice. \qed 

\begin{remark}\label{r:optcst}
For a bi-partite graph, the constant $2$ in
\eqref{e:majo} can be improved, in the sense of \eqref{e:reci0}, if
and only if $\Dc((\Delta_{\Er,   \theta})^{1/2})=
\Dc(d_G^{1/2}(Q))$.  
\end{remark} 

We finally provide an example. 

\begin{proposition}\label{p:bipart}
There is a simple bi-partite graph $K$ such that $\Dc(d^{1/2}_K(Q))\subsetneq
\Dc(\Delta^{1/2}_K)$. In particular, the constant $2$ in
\eqref{e:majo} is optimal in the sense of \eqref{e:reci0}. 
\end{proposition} 
\proof We start by constructing a complete bi-partite graph. Let
$K_{n,n}:=(\Er_n, \Vr_n)$ be defined as follows: $\Vr_n= \{1, \ldots
n\}\times \{1, \ldots n\}$ and such that  $\Er_n((k, i), (j,i))=0, \forall 
j,k \in \{1, \ldots n\}$ and $i=1,2$  and $\Er_n((k, 1), (j,2))=0, \forall
j,k \in \{1, \ldots n\}$. Consider now $f_n(x)=1$ on $K_{n,n}$. One has:
\begin{align}\label{e:fn}
\langle f_n, \Delta_K f_n\rangle =0, \quad \|f_n\|^2=2n, \quad
 \mbox{and} \quad  \langle f_n, d_K(Q) f_n\rangle= 2n^2.
\end{align}
Now take $K:=(\Er, \Vr)$, where $\Vr=\cup_{n\in \N\setminus\{0\}} \Vr_n$ and
$\Er(x,y):=\Er_n(x,y)$, if $x, y\in \Vr_n$, $\Er(x,y)=1$, if
$x=(1,2)\in \Vr_n$ and $y=(1,1)\in \Vr_{n+1}$, for all $n\geq 1$, and
$\Er(x,y)=0$ otherwise. 
\begin{align*}
\begin{array}{ccccccc}
\xymatrix{
&&&&&&{\bullet}\ar@{-}[r]\ar@{-}[dr]\ar@{-}[ddr]\ar@{-}[dddr]&{\bullet}\ar@{-}[dl]\ar@{-}[ddl]\ar@{-}[dddl]  &
\\
&&&&{\bullet}\ar@{-}[r]\ar@{-}[rd]\ar@{-}[ddr]&{\bullet}\ar@{-}[l]\ar@{-}[ld]\ar@{-}[ddl] 
&{\bullet}\ar@{-}[r]\ar@{-}[dr]\ar@{-}[ddr]&{\bullet}\ar@{-}[dl]\ar@{-}[ddl]&
\\
&&{\bullet}\ar@{-}[r]\ar@{-}[dr]&{\bullet}\ar@{-}[dl]
&{\bullet}\ar@{-}[r]\ar@{-}[dr]&{\bullet}\ar@{-}[dl]
&{\bullet}\ar@{-}[r]\ar@{-}[dr]&{\bullet} \ar@{-}[dl]&
\\
{\bullet}\ar@{-}[r]&{\bullet} \ar@{--}[r]&{\bullet}\ar@{-}[r]&{\bullet} \ar@{--}[r]
&{\bullet}\ar@{-}[r]&{\bullet}\ar@{--}[r]
&{\bullet}\ar@{-}[r]&{\bullet} \ar@{--}[r]&
} 
\\
\hspace*{-1.2cm} K_{1,1} \hspace*{1.7cm} K_{2,2} \hspace*{1.7cm}
K_{3,3}\hspace*{1.7cm} K_{4,4}
\\
\\
\mbox{ Graph of } K
\end{array}
\end{align*}
On the other hand, if $\Dc(d^{1/2}_K(Q))= \Dc(\Delta^{1/2}_K)$, the
uniform boundedness principle ensures that there are 
constants $a,b>0$, so that
\[ \langle f, d_K(Q) f\rangle \leq a \langle f, \Delta_K f\rangle + b\|f\|^2,
\quad \mbox{ for all } f\in \Cc_c(K).\]
This leads to a contradiction with \eqref{e:fn}, as $\langle f_n, \Delta_K
f_n\rangle=2$, for $n\geq 2$. Optimality follows by Proposition
and \ref{p:reci} and Remark \ref{r:optcst}.\qed

\section{Perturbation theory}\label{s:pertu}
We finally go into perturbation theory in order to obtain the stability
of the essential spectrum, of the inequality \eqref{e:pertutree}, and
of the asymptotic of eigenvalues. 

\begin{proposition}\label{p:pertutree}
Let $G=(\Vr, \Er, m, \theta)$ and $G_\circ=(\Vr,
\Er_\circ,m, \theta_\circ)$ be weighted graphs and $V:\Vr\to
\R$. Suppose that for all $\varepsilon>0$ there is $c_\varepsilon>0$
so that \eqref{e:tree1} holds true for $\Delta_{\Er_\circ,
  \theta_\circ}$. Suppose that there is
$\eta\in(0, 1)$ and $\kappa_\eta\geq0$, so that  
\begin{align}\label{e:pertutree}
|\langle f, V(Q)f \rangle|+ 2 \langle f, \Lambda(Q) f\rangle \leq \eta
\langle f, d_{G_\circ}(Q) f\rangle +\kappa_\eta\|f\|^2,
\end{align}
for all $f\in \Cc_c(\Vr)$, where
\begin{align}\label{e:lambda}
\Lambda(x):=  
\frac{1}{m^2(x)}\sum_{y\sim x} |\Er(x,y)-\tilde \Er(x,y)|.
\end{align}
Then, one has that:
\begin{enumerate}
\item The operator $\Hc:=(\Delta_{\Er,\theta}+V(Q))|_{\Cc_c(\Vr)}$ is
  bounded from below by some negative constant $-C$, in the form
  sense. We denote by  $\Hc_\Fr$ its Friedrichs extension. We have the
  equality of the form domains:
$\Dc(|\Hc_\Fr|^{1/2})= \Dc((\Delta_{\Er_\circ, \theta_\circ})^{1/2})$.
\item The three following assertions are equivalent:
\begin{enumerate}
\item[i)] The essential spectrum of
$\Hc_\Fr$ is empty,
\item[ii)] the essential spectrum of
$\Delta_{\Er_\circ, \theta_\circ}$ is empty,
\item[iii)] $\lim_{|x|\rightarrow   \infty} d_{G_\circ}(x)= +\infty$.
\end{enumerate}
\item Supposing that the essential spectrum of
$\Hc_\Fr$ is empty and that for all $\eta\in(0, 1)$
there is $\kappa_\eta\geq 0$, so that \eqref{e:pertutree} holds true,
then: 
\begin{align*}
\lim_{N \to
  \infty}\frac{\lambda_N(\Hc_\Fr)}{\lambda_N(d_{G_\circ}(Q))}=1.
\end{align*}
\end{enumerate}
\end{proposition}
 
\proof We start by noticing that, as in \eqref{e:majo},
\eqref{e:pertutree} ensures that:
\begin{align}\label{e:firstineq}
|\langle f, (\Delta_{\Er, \theta}+V(Q) - \Delta_{\Er_\circ,
  \theta_\circ}) f\rangle|& \leq  |\langle f, V(Q) f\rangle| + 2\langle
f, \Lambda(G) f\rangle\leq \eta \langle f,d_{G_\circ}(Q) f\rangle
+\kappa_\eta\|f\|^2 
\end{align}
for all $f\in \Cc_c(\Vr)$. Therefore, for all $\varepsilon>0$ and all
$\eta\in(0,1)$, satisfying \eqref{e:pertutree}, there is
$c_{\varepsilon, \eta}$ so that: 
\begin{align*}
(1- \eta -\varepsilon)\langle f, d_{G_0}(Q) f\rangle - c_{\varepsilon,
 \eta}\|f\|^2 &\leq 
\langle f, (\Delta_{\Er, \theta} +V(Q)) f\rangle 
\\
& \hspace{2cm}\leq
(1+ \eta+ \varepsilon) \langle f, d_{G_0}(Q) f\rangle+
c_{\varepsilon, \eta} \|f\|^2,
\end{align*}
for all $f\in \Cc_c(\Vr)$. This gives, directly, the first
point. Moreover, as above,  the second and third points follow by
Proposition \ref{p:compa}.  \qed

We now turn to the stability of the essential spectrum. 

\begin{proposition}\label{p:pertutree2}
Let $G=(\Vr, \Er, m, \theta)$ and $G_\circ=(\Vr,
\Er_\circ,m, \theta_\circ)$ be weighted graphs and $V:\Vr\to
\R$. Suppose that $\Dc(\Delta_{\Er_\circ, \theta_\circ}^{1/2})=
\Dc(d_{G_0}^{1/2}(Q))$ and that  
\begin{align}\label{e:pertutree1}
|V(x)|+\Lambda(x) = o(1+ d_{G_\circ}(x)), \mbox{ as } |x|\to \infty,
\end{align}
where $\Lambda$ is defined in \eqref{e:lambda}. Then
$\Hc:=\Delta_{\Er,\theta}+V(Q)|_{\Cc_c(\Vr)}$ is bounded from
below by some negative constant $-C$, in the form sense. We denote by 
$\Hc_\Fr$ its Friedrichs extension. We obtain that  
$\Dc(|\Hc_\Fr|^{1/2})=
\Dc(d_{G_0}^{1/2}(Q))$
and 
$\sigma_{\rm ess}(\Delta_{\Er_\circ, \theta_\circ})= \sigma_{\rm
  ess}(\Hc_\Fr)$. 
\end{proposition} 

\proof The uniform boundedness theorem and \eqref{e:pertutree1}
implies that $\Hc$ is bounded from below by some
negative constant $-C$, in the form sense. We consider its Friedrichs
extension.  Next, KLMN's Theorem, e.g., \cite[Theorem X.17]{RS}
ensures that $\Dc(|\Hc_\Fr|^{1/2})=
\Dc(d_G^{1/2}(Q))$.  

By Weyl's Theorem, e.g., \cite[Theorem
XII.1]{RS}, it is enough to show 
that the difference  of the resolvents is compact. As we have to work
with forms,  one should be careful with the resolvent equation. We
give a complete proof and refer to \cite{GG} for more discussions of
this matter. To lighten notation, we set $H_0:=\Delta_{\Er_\circ,
  \theta_\circ}$ and $H:=\Delta_{\Er,\theta}+V$. To start off, we give
a    rigorous meaning to   
\begin{equation}\label{33}
(H+\rmi)^{-1} - (H_0+\rmi)^{-1}=
(H+\rmi)^{-1}
(H_0 - H)
(H_0+\rmi)^{-1}.
\end{equation}
Since $\Gr:=\Dc((H+C)^{1/2})= \Dc((H_0)^{1/2})$, both operators
extend to an element of $\Bc(\Gr, \Gr^*)$. Here we use the Riesz
lemma to identify $\Hr$ with its anti-dual $\Hr^*$. We denote these
extensions with a tilde.    

We have $(H_0+\rmi)^{-1*}\Hr\subset\Gr$. This allows one to
deduce that ($H_0+\rmi)^{-1}$ extends to a unique continuous operator
$\Gr^*\rightarrow\Hr$. We denote it for the moment by $R$. From
$R(H_0+\rmi)u=u$ for $u\in \Dc(H_0)$ we get, by density of 
$\Dc(H_0)$ in $\Gr$ and continuity, $R(\widetilde{H_0}+\rmi)u=u$ for
$u\in\Gr$. In particular
\begin{equation*}
(H+\rmi)^{-1}= R(\widetilde{H_0}+\rmi)(H+\rmi)^{-1}.
\end{equation*}
Clearly,
\begin{equation*}
(H_0+\rmi)^{-1}= (H_0+\rmi)^{-1}(H+\rmi)(H+\rmi)^{-1}=
R(\widetilde{H}+\rmi)(H+\rmi)^{-1}.
\end{equation*}
We subtract the last two relations to obtain that
\begin{equation*}
(H_0+\rmi)^{-1}-(H+\rmi)^{-1}=R(\widetilde{H}
-\widetilde{H_0})(H+\rmi)^{-1}. 
\end{equation*}
Since $R$ is uniquely determined as the extension of $(H_0+\rmi)^{-1}$ to a
continuous map $\Gr^*\rightarrow\Hr$, one may keep the notation
$(H_0+\rmi)^{-1}$ for it. With this convention, the rigorous version of
(\ref{33}) that we use is:
\begin{equation}\label{34}
(H_0+\rmi)^{-1}-(H+\rmi)^{-1}=(H_0+\rmi)^{-1}(\widetilde
H-\widetilde H_0)(H+\rmi)^{-1}. 
\end{equation}
Therefore, to prove the equality of the essential spectra, it is
enough to show that $\widetilde 
H-\widetilde H_0$ is a compact operator from $\Gr$ to $\Gr^*$. 
By \eqref{e:firstineq}, one gets:
\begin{align*}
|\langle(1+d_{G_0})^{-1/2}(Q) f, (\Delta_{\Er, \theta}+V(Q)- \Delta_{\Er_\circ,
  \theta_\circ}) (1+d_{G_0})^{-1/2}(Q) f\rangle| &\leq
\\  &\hspace*{-3cm}\langle
(1+d_{G_0})^{-1/2}(Q) f, (V+2\Lambda)(Q) (1+d_{G_0})^{-1/2}(Q)f\rangle,
\end{align*}
for all $f\in \Cc_c(\Vr)$. Now, by hypothesis \eqref{e:pertutree1},
$(V+2\Lambda)(1+d_{G_0})^{-1}(Q)$ is compact in $\ell^2(\Vr, m^2)$. We
conclude by using  Proposition \ref{p:compact}.\qed

\begin{remark}
Note that \eqref{e:pertutree} and \eqref{e:pertutree1} allow us to
consider some potentials $V(Q)$ that are unbounded from below,
whereas $\Delta_{\Er,    \theta}+V(Q)$ is bounded from below. This
is due to the fact that we know the form-domain explicitly.
\end{remark} 

\appendix
\section{The $\Cc^1-$Regularity} \label{a:regu}
\setcounter{equation}{0}  
We start with some generalities. Given a bounded operator $B$ and a
self-adjoint operator $A$ acting in a Hilbert space $\Hr$, one says
that $B\in \Cc^k(A)$ if $t\mapsto e^{-\rmi tA}B e^{\rmi tA}$ is strongly
$\Cc^k$. Given a closed and densely defined operator 
$B$, one says that $B\in \Cc^k(A)$ if for some (hence any) $z\notin
\sigma(B)$, $t\mapsto e^{-\rmi tA}(B-z)^{-1} e^{\rmi tA}$ is strongly
$\Cc^k$. The two definitions 
coincide in the case of a bounded self-adjoint operator. We recall a
result following from Lemma 6.2.9 
and Theorem 6.2.10 of \cite{ABG}. 
\begin{theorem}\label{th:abg} 
Let $A$ and $B$ be two self-adjoint operators in the Hilbert space
$\Hr$. For $z\notin \sigma(A)$, set $R(z):=(B-z)^{-1}$. The following
points are equivalent to $B\in\Cc^1(A)$:  
\begin{enumerate} 
\item For one (then for all) $z\notin \sigma(B)$, there is a finite
$c$ such that 
\begin{align}\label{e:C1a} 
|\langle A f, R(z) f\rangle - \langle R(\overline{z}) f, Af\rangle| \leq c 
\|f\|^2, \mbox{ for all $f\in\Dc(A)$}. 
\end{align} 
\item 
\begin{enumerate} 
\item [i)]There is a finite $c$ such that for all $f\in \Dc(A)\cap\Dc(B)$: 
\begin{equation}\label{e:C1b} 
|\langle Af, B f\rangle- \langle B f, Af\rangle|\leq \,
 c\big(\|B f\|^2+\|f\|^2\big). 
\end{equation} 
\item [ii)] For some (then for all) $z\notin \sigma(B)$, the set
$\{f\in\Dc(A),  R(z)f\in\Dc(A)$ and $R(\overline{z})f\in\Dc(A)
\}$ is a core for $A$. 
\end{enumerate} 
\end{enumerate} 
\end{theorem} 
Note that the condition ii) could be delicate to check (see
\cite{GeorgescuGerard}). We mention \cite{GoleniaMoroianu}[Lemma A.2]
and \cite{GerardLaba}[Lemma 3.2.2] to overcome this subtlety. 

Note that \eqref{e:C1a} yields that the commutator $[A, R(z)]$ extends to a
bounded operator, in the form sense. We shall denote the extension by
$[A, R(z)]_\circ$. In the same way, since $\Dc(B)\cap \Dc(A)$ is dense
in $\Dc(B)$,  \eqref{e:C1b} ensures that the commutator $[B, A]$
extends to a unique element of $\Bc\big(\Dc(B), \Dc(B)^*\big)$ denoted
by $[B,   A]_\circ$. Moreover, when $B\in \Cc^1(A)$, one has: 
\begin{eqnarray}\label{e:facto}
\big[A, (B-z)^{-1}\big]_\circ =\quad  \underbrace{(B-z)^{-1}}_{\Hr
  \leftarrow \Dc(B)^*}\quad  \underbrace{[B, A]_\circ}_{\Dc(B)^*\leftarrow
  \Dc(B)} \quad \underbrace{(B-z)^{-1}}_{\Dc(B)\leftarrow \Hr}.
\end{eqnarray} 
Here we use the Riesz lemma to identify $\Hr$ with its anti-dual
$\Hr^*$. 

It turns out that an easier characterization is available if 
the domain of $B$ is conserved under the action of the $C_0-$group 
generated by $A$. 
\begin{theorem}(\cite[p.\ 258]{ABG})\label{th:abg2} 
Let $A$ and $B$ be two self-adjoint operators in the Hilbert space
$\Hr$ such that $e^{itA}\Dc (B)\subset\Dc (B)$, for all $t\in\R$. 
Then, $B\in\Cc^1(A)$ if and only if 
\[\liminf_{t\rightarrow 0^+} \left\|[B, e^{\rmi t
    B}/t]\right\|_{\Bc(\Dc (B), \Dc (B)^*)}<\infty.\] 
\end{theorem}
Note that $e^{itA}\Dc (B)^*\subset\Dc (B)^*$ by duality. 

\section{The Helffer-Sj\"{o}strand's formula} \label{a:comm}
\setcounter{equation}{0}  

We present briefly  the  Helffer-Sj\"{o}strand's formula. We refer to
\cite{GoleniaJecko}[Appendix B] and \cite{BG}[Appendix A] (see also
\cite{DerezinskiGerard, HunzikerSigalSoffer, M}) for commutator
expansion. We first recall some well-known facts about almost analytic
extensions. For $\rho\in\R$,   let $\Sc^\rho$ be the class of function
$\varphi\in\Cc^\infty(\R;\C)$ such that       
\begin{eqnarray}\label{eq:regu} 
\forall k\in\N, \quad C_k(\varphi) :=\sup _{t\in\R}\, \langle
t\rangle^{-\rho+k}|\varphi^{(k)}(t)|<\infty.
\end{eqnarray} 
Equipped with the semi-norms defined by (\ref{eq:regu}), $\Sc^\rho$ is
a Fr\'echet space. Leibniz' formula implies the continuous embedding:
$\Sc^\rho\cdot  \Sc^{\rho'} \subset  \Sc^{\rho+\rho'}$.  
We shall use the following result, e.g., \cite[Appendix C.2]{DerezinskiGerard}. 
 
\begin{lemma}\label{l:dg} 
Let $\varphi\in\Sc^\rho$ with $\rho\in\R$. For all   
$l\in \N$, there is a smooth function  $\varphi^\C:\C \rightarrow \C$,
such that: 
\begin{eqnarray} 
\label{eq:dg1} \varphi^\C|_{\R}=\varphi,\quad &&\left|\frac{\partial 
  \varphi^\C}{\partial  \overline{z}}(z) \right|\leq c_1 \langle \Re(z) 
  \rangle^{\rho-1 -l} |\Im(z)|^l\\\label{eq:dg2} 
&& \supp\, \varphi^\C \subset\{x+\rmi y,  |y|\leq c_2 \langle 
  x\rangle\},\\\label{eq:dg3} && \varphi^\C(x+\rmi y)= 0, \mbox{ if } 
  x\not\in\supp\,\varphi .  
\end{eqnarray}  
for some constants $c_1$, $c_2$ depending on the semi-norms 
\eqref{eq:regu} of $\varphi$ in $\Sc^\rho$ and not on $\varphi$.   
\end{lemma}  
One calls $\varphi^\C$ an \emph{almost analytic extension} of $\varphi$.
Let $A$ be a self-adjoint operator, $\rho < 0$ and $\varphi\in
\Sc^{\rho}$. By functional calculus, one has $\varphi(A)$
bounded. The Helffer-Sj\"{o}strand's formula, e.g., 
\cite{HelfferSjostrand} and \cite{DerezinskiGerard},
gives that for all almost analytic extension of $\varphi\in\Sc^\rho$,
with $\rho<0$, we have:  
\begin{eqnarray}\label{eq:int} 
\varphi(A) = \frac{\rmi}{2\pi}\int_\C\frac{\partial \varphi^\C }{\partial 
\overline{z}}(z-A)^{-1}dz\wedge d\overline{z}. 
\end{eqnarray}  
Note that the integral exists in the norm topology, by \eqref{eq:dg1} 
with $l=1$ and by taking in account the domain of integration given in
\eqref{eq:dg2}.

\end{document}